\documentclass[10pt]{article}
\usepackage{amsmath,amsfonts,amsthm,amssymb,amscd}

\catcode`\é=\active\def é{\'e} \catcode`\è=\active\def è{\`e}
\catcode`\ç=\active\def ç{\c{c}} \catcode`\à=\active\def à{\`a}
\catcode`\ê=\active\def ê{\^e} \catcode`\ù=\active\def ù{\`u}
\catcode`\ô=\active\def ô{\^o} \catcode`\û=\active\def û{\^u}
\catcode`\î=\active\def î{\^{\i}} \catcode`\â=\active\def â{\^a}
\catcode`\ö=\active\def ö{\"o}

\newcommand{\be}{\begin{equation}}
\newcommand{\ee}{\end{equation}}
\newcommand{\bs}{\begin{split}}
\newcommand{\es}{\end{split}}
\newcommand{\ba}{\begin{align}}
\newcommand{\ea}{\end{align}}
\newcommand{\basl}[1]{\begin{align}\begin{split}\label{#1}}
\newcommand{\bas}{\begin{align}\begin{split}}

\newtheorem{theo}{Theorem}[section]
\newtheorem{prop}[theo]{Proposition}
\newtheorem{lemm}[theo]{Lemma}

\newcommand\R{\mathbb{R}}
\newcommand\C{\mathbb{C}}

\title{Hamiltonian systems and semiclassical dynamics \\for interacting spins in QED}

\author{L. Amour and J. Nourrigat}

\date{Université de Reims, France}

\begin{document}

\maketitle

\begin{abstract}
\noindent
In this article, we consider fixed spin$-1/2$ particles interacting through the quantized electromagnetic field in a constant magnetic field. We give approximate evolutions of coherent states. This uses spins-photon classical Hamiltonian mechanics. These approximations enable to derive that the approximate  average fields and spins  follow Maxwell-Bloch equations with a current density coming from spins. In addition, we obain a law concerning the evolution of the approximate average number of photons. Next, we provide stationary points of the spins-photon Hamiltonian when the spin particles belong to an orthogonal plane to the constant magnetic field. This allows for the construction of quasimodes with spins colinear to the constant magnetic field. Finally, a quasimode with an arbitrary high accuracy is built up and its first order radiative correction is computed.
\end{abstract}

\parindent=0pt

\


\


\tableofcontents

\parindent = 0 cm

\parskip 10pt
\baselineskip 17pt

\section{The model.}\label{s1}

The aim of this work is to initiate a semiclassical study for a Hamiltonian operator modelling the interaction between quantized electromagnetic field and $N$ fixed spin-$1/2$ particles in a constant magnetic field. To achieve that goal, a Hamiltonian system is defined, playing the same role for this Hamiltonian operator as the usual system for the Schr\"odinger operator. This Hamiltonian is first recalled below in  (\ref{1.14}) and (\ref{1.15})  (see Reuse
\cite{REU}, H\"ubner-Spohn \cite{H-S}, Derezi\'nski-Gérard \cite{D-G}  Sections 3.2-3.3).
The Hilbert space associated with this Hamiltonian is the completed tensor product ${\cal H}_{ph} \otimes {\cal H}_{sp}$.

The Hilbert space  ${\cal H}_{ph}$ for photons may be viewed as the symmetrized Fock space ${\cal F}_s (H _{\bf C} )$ associated with the complexified of some real Hilbert space $H$ defined by Lieb-Loss \cite{L-L}. This space $H$  is the space of mappings
 $f = (f_1, f_2, f_3)$ from $\R^3$ to $\R^3$
with $f_j$ belonging in $L^2({\R}^3)$, taking real values and satisfying,
\be\label{1.1}
k_1f_1(k) + k_2f_2(k) + k_3f_3(k) = 0\quad {\rm a.e.}\ee
This space is equipped with the norm,
\be\label{1.2} |f|^2 = \sum_{j=1}^3 \int_{\R^3} |f_j(k)|^2dk. \ee

The Fock space ${\cal F}_s (H _{\bf C} )$  definition is reminded  in Section \ref{s6}.
For each $h>0$, the space ${\cal H}_{ph}$ may also be viewed as a
$L^2(B , \mu _{B , h/2})$ space where $B$ is some Banach space related to $H$ and containing $H$,
and $\mu _{B , h/2}$ is some gaussian measure with variance $h/2$ on the Borel $\sigma-$ algebra of $B$.
The hypotheses that  should  be fulfilled by  $B$ are recalled in \cite{A-J-N} or \cite{A-L-N-2}. When these assumptions are satisfied then $(i, H, B)$ is called Wiener space where $i$ the injection from $H$ into $B$.

The space ${\cal H}_{sp}$ for the spin particles is denoted by $( \C^2  )^{\otimes N}$.

In the space  ${\cal H}_{ph}$, the definition of the model involves three kinds of operators: the number operator $N$, the free photons  energy operator $H_{ph}$ and operators at each point $x\in \R^3$ associated with the three components of the magnetic field. These operators are denoted by $B_m (x)$, $1 \leq m \leq 3$, and for the electric field, it is denoted by  $E_m (x)$, $1 \leq m \leq 3$.
Each of these operators is depending on the semiclassical parameter $h>0$ which is sometimes not explicitly written.

Within the Fock space formalism, the number operator $N$ and the free photons Hamiltonian
 $H_{ph}$  are defined by,
%---
\be\label{1.3}  N= {\rm d}\Gamma (I),\qquad H_{ph} = h {\rm d}\Gamma (M), \ee
%---
$M$ being the multiplication operator by $\omega (k) = |k|$ with domain $D(M) \subset H$, ${\rm d}\Gamma $ is the standard operator (see \cite{RE-SI} or (\ref{6.17}) below) and $h>0$ is the semiclassical parameter. These equalities classically define selfadjoint operators (see \cite{RE-SI}).

Using  Wiener spaces formalism, the Wick symbol  (see \cite{A-J-N}\cite{A-L-N-1}\cite{A-L-N-2}) of the number operator $N$ is $(q, p) \mapsto ( |q|^2 + |p|^2)
/2h$. The Wick symbol of the operator $H_{ph}$ is
%---
\be\label{1.4} H_{ph} (q , p) = {1\over 2} \int _{\R^3} |k| \Big [ |q(k)|^2 +  |p(k)|^2 \Big ] dk. \ee
%----

With the notations of \cite{A-J-N}, the operators $B_m (x)$ and $E_m (x)$ are associated with functions on  $H^2$, denoted by $(q , p) \mapsto B_m (x, q, p)$ and
$(q , p) \mapsto E_m (x, q, p)$, being the Weyl and Wick symbols of these operators. These symbols are continuous linear forms on $H^2$ and are denoted by,
%-----
\be\label{1.5} B_m (x, q, p) =  a_m(x) \cdot q + b_m(x) \cdot p,\qquad x\in \R^3, \ \ \ 1 \leq m \leq 3,
 \ \ \ (q , p)\in H^2,  \ee
%----
with the scalar product of  $H$. Here $a_m(x)$ and $b_m(x)$ are elements of $H$. Thus, $a_m(x)$ and $b_m(x)$ are mappings from $\R^3$ into itself and are defined by,
 %----
\be\label{1.6} a_m (x) (k) =
 {\chi(|k|)|k|^{1\over 2} \over (2\pi)^{3\over 2}} \sin ( k \cdot x )
{k\wedge e_m \over |k|}\ee
%---
\be\label{1.7} b_m (x)   (k) =
 {\chi(|k|)|k|^{1\over 2} \over (2\pi)^{3\over 2}} \cos  ( k \cdot x )
{k\wedge e_m \over |k|}, \ee
%---
where $\chi $ is a function belonging to ${\cal S} (\R)$  vanishing in a neighborhood of the origin and $(e_1, e_2, e_3)$
is the canonical basis of  $\R^3$. The fact that  $\chi $ vanishes near  $0$ is used to apply results of \cite{A-L-N-2} only involved in Lemma \ref{l3.6} and it is also used in the  construction of a quasimode in Section \ref{s6}. The link between the symbol $B_m (x, p, q)$ and the operator $B_m(x)$ may be defined either by the Weyl calculus studied in \cite{A-J-N} and \cite{A-L-N-1} or with the Segal field formalism.

In \cite{A-L-N-1}, an operator denoted by $Op_h^{weyl} (F)$ is associated with any continuous linear form $F$ on $H^2$, for any $h>0$ which, in its initial definition, maps some space
${\cal D}$ into itself and is thereafter extended to a bounded operator from the space $W_1$
defined in Section \ref{s3} into ${\cal H}_{ph}$. When setting
%---
\be\label{1.8} F_{a , b} (q , p) = a \cdot q + b \cdot p, \ee
%----
the operator $Op_h ^{weyl} (F) $ is usually called Segal field and is denoted by $\Phi_S(a+ib)$
(see \cite{RE-SI}). The point in \cite{A-J-N} and \cite{A-L-N-1} is to be able to define operators $Op_h^{weyl} (F)$ for others functions $F$ which are not linear, but this is only used in Lemma \ref{l3.6}.
With these  notations, one has,
%----
\be\label{1.9} B_m (x) = Op_ h^{weyl} (F_{ a_m (x), b_m (x)} ) = \Phi_S (a_m (x) +i b_m (x)).\ee
%-------

In the following, the operator is denoted by  $B_m (x)$ and its symbol is denoted by $B_m (x , \,\cdot\,)$. Concerning the symbol $E_m (x, q, p)$ of the operator $E_m(x)$, it is related to $B_m (x, q, p)$ through the helicity operator $J$ mapping $H^2$ into $H^2$ and defined by,
%----
\be\label{1.10} J(q, p) (k) = \left ( {k\wedge q(k) \over |k| } , {k\wedge p(k) \over |k| } \right ), \qquad
k\in \R^3 \setminus \{ 0 \}.\ee
%-----
The symbol $  E_j (x, q, p)$ is defined by,
%----
\be\label{1.11}  E_j (x, q, p) =  - B_j (x , J(q, p)).\ee
%---

Operators in ${\cal H}_{sp}$ use in particular Pauli
 matrices $\sigma_j$ ($1\leq j \leq 3$),
%----
\be\label{1.12} \sigma_1 = \begin{pmatrix}  0 & 1 \\ 1 & 0   \end{pmatrix},
\qquad
   \sigma_2 = \begin{pmatrix} 0 & -i \\ i & 0   \end{pmatrix},
 \qquad
   \sigma_3 = \begin{pmatrix}  1 & 0 \\ 0 & -1  \end{pmatrix}.\ee
%----
For all $\lambda \leq N$ and for any $j\leq 3$,   $\sigma_j^{[\lambda]}$
denotes the following operator in  ${\cal H}_{sp}$,
%---
\be\label{1.13} \sigma_j^{[\lambda]} = I \otimes \cdots \sigma_j \cdots \otimes I,\ee
%---
where $\sigma_j$ is located at the $\lambda ^{th}$ position.

We assume that there are $N$ fixed spin$-1/2$ particles at points $x_{\lambda}$ in $\R^3$
($1\leq \lambda \leq N$).
Denoting the constant magnetic field by $\beta = (\beta_1, \beta_2, \beta _3)$,
the system constituted with these particles and the quantized magnetic field is governed by the operator in ${\cal H}_{ph} \otimes {\cal H}_{sp}$ defined by,
%----
\be\label{1.14} H (h) = H_0 + h H_{int},\qquad H_0 = H_{ph}\otimes I ,\ee
%----
where
%----
\be\label{1.15} H_{int} = \sum _{\lambda = 1}^N \sum_{j=1}^3
( \beta _j + B_j(x_{\lambda}) ) \otimes  \sigma _j ^{[\lambda]}.\ee
%---
It is recalled in \cite{A-L-N-2} (Section 4) that it defines a selfadjoint operator with domain
 $D(H_{ph}) \otimes {\cal H}_{sp}$.

In order to describe the same interaction in the semiclassical limit, a symplectic manifold is introduced in
Section \ref{s2}, namely $H^2 \times \Omega $, where
$H$ is the space of Lieb-Loss and $\Omega$ is an orbit of the coadjoint  representation of the group
 $SU(2^N)$. On this symplectic manifold, we introduce a Hamiltonian system  modelling the interaction  between  photons and spins in the semiclassical regime.
Let us mention a similarity with the Hamiltonian system of  Bolte-Glaser \cite{B-G} modelling the evolution of a particle whose spin interacts with a given magnetic field (semiclassical limit of Pauli equations). In Section \ref{s2}, an existence and uniqueness result for weak and strong solutions to the Hamiltonian system is given.

A linear form $\pi _a$ on the Lie algebra ${\cal G}$ of the Lie group $SU(2^N)$ is associated with each element $a$ in ${\cal H}_{sp}$
defined by,
%---
\be\label{1.16} <\pi_a ,  Y > =i^{-1}  < Ya, a>,\qquad Y \in {\cal G}.\ee
%---
The unit sphere of ${\cal H}_{sp}$ is thus a circle bundle  above the coadjoint orbit $\Omega$.

In (\ref{2.15}), with each solution $(q(t), p(t) , \omega (t))$ of the Hamiltonian system and with any $a_0$ such that
$\pi _{a_0} = \omega (0)$, it is associated a flow $a(t)$  in  ${\cal H} _{sp}$ satisfying $\pi _{a(t)} = \omega (t)$,
for all $t$.

The first application of this Hamiltonian system is the study of the approximation as $h$ goes to $0$
of  evolutions of coherent states. Coherent states may be very simply defined as elements of ${\cal F}_s (H _{\bf C} )$.
For all $X=(a , b)$ in $H^2$ and for any $h>0$,  $\Psi_{X ,h}$ denotes the corresponding coherent state, element of  ${\cal F}_s (H_{\bf C})$ defined by,
%----
\be\label{1.17} \Psi_{(a , b) , h}  = \sum _{n\geq 0}
 {e^{-{|a|^2+ |b|^2  \over 4h}} \over (2h)^{n/2} \sqrt {n!} } (a+ib) \otimes \cdots \otimes (a+ib).\ee
%----
If ${\cal H } _{ph } $ and $ L^2(B , \mu_{B , h/2})$ are identified then the coherent state
$\Psi_{(a , b) , h}$ is explicitly written in (\ref{3.5}) below.

Let us fix $ X_0 = (q_0 , p_0)$ in $H^2$ and  $a_0$ in the unit sphere of ${\cal H} _{sp}$.
Our goal is to give an approximation when $h$ goes to $0$ of
%---
\be\label{1.18} u_h (t) =   e^{ - {it\over h} H (h) } \Big ( \Psi_{X_0 , h} \otimes a_0 \Big ).\ee
%---

\begin{theo}\label{t1.1} Let $H(h)$ be the Hamiltonian defined in (\ref{1.14})(\ref{1.15}).
Let $X_0 = (q_0 , p_0)$ belongs to $H^2$ and $a_0$ in ${\cal H}_{sp}$ with a unit norm.
Let $(q(t) , p(t) , \omega (t) )
= (X(t) , \omega (t) ) $  be the weak solution to the Hamiltonian system (\ref{2.6})(\ref{2.7})(\ref{2.10}) with initial data $(q_0 , p_0, \pi_{a_0})$.
Let   $a(t)$ be the solution to the differential system (\ref{2.15}) with initial data $a_0$.
Set,
%---
\be\label{1.19} \varphi (t) = \int_0^t \Big [ {1\over 2} \Big (
   p(s) \cdot q'(s)- q(s) \cdot p'(s)  \Big )-  H( q(s), p(s), \omega (s) )  \Big ] ds,\ee
%---
where $H(q, p, \omega)$  is the Hamiltonian function defined in (\ref{2.5}).
Let $u_h(t)$ be defined in  (\ref{1.18}) and
%---
 \be\label{1.20}v_h(t) = e^{ {i\over h} \varphi (t)}  \Big ( \Psi_{X(t) , h} \otimes a (t) \Big ).\ee
%---
Then
%---
\be\label{1.21}  \Vert u_h(t) - v_h(t) \Vert  \leq C(t) h^{1/2} .\ee
%---
The above norm is the one of ${\cal H}_{ph}\otimes {\cal H}_{sp}$ and the above constant $C(t)$ is bounded on every compact set of $\R$.
\end{theo}

This result is proved in Section \ref{s3} also using a Sobolev space norm being necessary for applications.

One may compare these results to those in finite dimension of  Combescure Robert \cite{C-R}  (also see
\cite{RO} and Hagedorn \cite{HA}).

There is here a different approach which is somewhat close to the one of Bolte-Glaser \cite{B-G}. Even if the spin terms are not part of what one may consider as the principal symbol, these are taken into account in the Hamiltonian system. This enables the Hamiltonian to model the interaction between photons and matter and not only the free evolution of the photons.

We are next interested with the time evolution of observables. For example, with the initial data of Theorem \ref{t1.1}, we consider the average number of photons at time $t$ given by $< N u_h(t), u_h (t)>$. One of the points of
Theorem \ref{t1.1} (adapted with Sobolev spaces) is to approximate the average number of photons  by $<N v_h(t), v_h (t)>$ with an error of ${\cal O} (h^{1/2})$.
This also holds true for the three components of the electric field and the magnetic field at each point $x\in \R^3$ and at any time $t$. This remains still valid for the three components of the spin (viewed as a vector in $\R^3$) of each
particle at time $t$. The point of these approximations is that, these average values of these observables taken on coherent states can be explicitly  computed. These computations enable to describe approximate laws for the evolution of the above observables.
For the case of the magnetic field, one finds Maxwell equations with a divergence free current, without charge, modelling the spin (see (\ref{5.3})-(\ref{5.6})). For the case of spins viewed as vectors in  $\R^3$, one recovers
the equations of F. Bloch \cite{BL} (1946) (see \ref{5.9}). For the evolution of the average number of photons, one finds a time evolution law, may be new,  that one may express as follows if there is a unique spin particle. If the photons of the coherent states are in a given circularly polarized state then then the derivative of the the average number of photons at time $t$ is equal to the scalar product of the spin at this time with the electric field at the point where the particle is located, multiplied  by $+1$ or $-1$ according to the circular polarization direction (see (\ref{5.10})). It is clearly an approximate law, but one may expect to give more precise laws with quantum corrections in a next article. These points are detailed in Sections \ref{s4} and \ref{s5}.

Next, we consider fixed points of the Hamiltonian system. We show that, if there are  $N$ spin$-1/2$ particles then there are $2^N$ fixed points of the Hamiltonian system provided that all the particles are located in the same plane orthogonal to the constant magnetic field. Energy levels are given in (\ref{6.5})(\ref{6.6}).
In order to give them a more physical meaning, the function
$\chi$ in (\ref{1.6})(\ref{1.7}) has to tend to $1$, which is possible since the hypothesis on $\chi $ vanishing at the origin is useless here. Taking the limit
$\chi $ goes to $1$, while subtracting a common term going to infinity to all  energy levels, one then finds for critical energy levels,
%---
\be\label{1.22} E =h |\beta | \sum _{\lambda = 1}^N \varepsilon _{\lambda }
 + {h^2 \over 8 \pi } \sum _{\lambda \not=  \mu  } { \varepsilon_{\lambda } \varepsilon_{\mu }
\over |x_{\lambda } - x_{\mu  }|^3 },\ee
%---
where $\varepsilon_{\lambda } = \pm 1$ depending if the spin $\lambda $ of the considered critical point is orientated in one direction or the other. We still suppose that the
points $x_{\lambda }$  are in the same plane orthogonal to the constant field. One recognizes the Ising model.

Finally, we construct a quasimode, which may be have a connection with the ground state whose existence is proved in \cite{H-S} or \cite{D-G} (for massive photons) and \cite{G} (for massless photons).
See \cite{PA-U} for a construction of  quasimodes with coherent states.
The quasimode is constructed in Section \ref{s6} with the Fock
formalism under the form, %---
$$ u_h \sim  u_0 + h^{1/2} u_{1} + h u_2 + \cdots $$
%---
where the $u_j$ are elements of ${\cal H}_{ph } \otimes {\cal H}_{sp }$ and independent on $h$ (if considered as elements  of the Fock space). There are no hypotheses on the position of the particles for this  quasimode. For all  $x$ in $\R^3$, one has, assuming that the constant magnetic filed is
$\beta=(0,0,|\beta|)$,
%----
\be\label{1.23} < (B_m (x)\otimes I)  ( u_0 + h^{1/2} u_{1} ), ( u_0 + h^{1/2} u_{1} ) >  = \ee
%-----
$$ h  (2 \pi )^{-3} \sum _{\lambda = 1}^N \int _{\R^3} |\chi (k)|^2
\ \cos (k \cdot ( x - x_{\lambda } )) \ { ( k \wedge e_m ) \cdot ( k \wedge e_3 )\over |k|^2 } dk.   $$
%----
The corresponding electric field is zero.  Denoting by ${\bf B } (x)$ the vector with the above three components, it is seen that,
%---
\be\label{1.24} {\rm div } {\bf B } = 0,\qquad{\rm rot } {\bf B } = -  h e_3 \wedge {\rm grad} \Phi (x),\ee
%---
\be\label{1.25}  \Phi (x) = \sum _{\lambda = 1}^N \rho (x - x_{\lambda }), \qquad
\rho (x) = (2 \pi )^{-3}  \int _{\R^3} |\chi (k)|^2  \ \cos (k \cdot  x) dk.\ee
%---
One notes the role of the function
 $\chi$: the spin
is not exactly a point-like particle. One recovers the stationary Maxwell equations in the first order approximation
of the magnetic field associated with the quasimode. The approximate eigenvalue  has an asymptotic expansion $\lambda (h)  \sim \lambda_1 h + \lambda_2 h^2 + \dots$ and we have explicit formulas for the first two terms.

\section{Photons and spins Hamiltonian mechanics.}\label{s2}

Let us begin by defining the symplectic manifold of our system.

Let $G = SU (2^N) $ be the group of unitary matrices in ${\cal H}_{sp}$.
The  Lie algebra ${\cal G} $ of this group is the algebra of traceless antihermitian maps.
For all $a$ in the unit sphere of ${\cal H}_{sp}$, let $\pi_a$ be the linear form on  ${\cal G}$  defined in (\ref{1.16}).
Let $\Omega \subset {\cal G} ^{\star} $ be the set of all  linear forms $f$
written as $\pi_a$ with $a$ in the unit sphere of  ${\cal H}_{sp}$.
When $a$ is any element of this sphere, one notes that $\Omega $
is the orbit $\pi_a$ of the coadjoint representation,
that is to say, the set of  linear forms on ${\cal G}$ written as,
%----
\be\label{2.1} {\cal G} \ni Y \mapsto i^{-1} < Y ga , ga >,\ee
%---
with  $g\in G$.  It is known that $\Omega $ is a  symplectic manifold. For all $Y$ in $ {\cal G}$, let
$\varphi _Y$ be the function on $\Omega$ defined by  $\varphi_Y (\omega) = < \omega , Y> $. Then,
the Poisson bracket of $\varphi_Y$ and $\varphi_Z$ ($Y$ and $Z$ in ${\cal G} $) satisfies,
%---
\be\label{2.2} \{ \varphi _Y , \varphi _Z \} = h^{-1} \varphi _{[Y, Z]}.\ee
%------

 When $H$ is the Hilbert space of Lieb-Loss defined in Section 1, the space $H^2$ is equipped
with the symplectic form defined by,
%--
\be\label{2.3} \sigma ( ( q, p) , (q' , p')) = ( q' , p) - (q , p').\ee
%----
The  symplectic manifold of our system is,
%----
\be\label{2.4} V = H^2 \times \Omega.\ee
%-----

The Hamiltonian function $ H_{ph} (q , p)$ for photons  used here is defined in (\ref{1.4}).
The functions $B_m (x, q, p)$ and the operators
$\sigma _m^{[\lambda ]}$ are defined in Section 1.
We denote by $\beta = ( \beta _1 , \beta _2 , \beta _3)$ the constant magnetic field.

The Hamiltonian function for photons in interaction with  $N$ spin$-1/2$ particles fixed at
the points $x_{\lambda}$ ($1 \leq \lambda \leq N$) is,
%---
\be\label{2.5} H(q, p, \omega ) = H_{ph} (q , p) +h  \sum _{\lambda = 1}^N \sum _{m=1}^3
( \beta _m + B_m (x_{\lambda }, q, p) ) <\omega , i \sigma _m^{[\lambda ]}  >,\qquad (q, p)\in H^2, \ \ \ \ \omega  \in \Omega.\ee
%----
The above bracket denotes the duality between ${\cal G}'$ and ${\cal G}$.

Consequently, a solution $(q(t), p(t), \omega (t))$ to the  Hamiltonian system
should satisfy,
%----
\be\label{2.6} {d q(t) \over dt }  =  M p (t) +h  \sum _{\lambda = 1}^N \sum _{m=1}^3
 b_m (x_{\lambda } )   \ < \omega (t)) , i\sigma _m^{[\lambda ]}>\ee
 %------
 %----
\be\label{2.7} {d p(t) \over dt }  = - M q(t) - h \sum _{\lambda = 1}^N \sum _{m=1}^3
 a_m (x_{\lambda }  )   \ < \omega (t)) , i \sigma _m^{[\lambda ]}>.\ee
 %------
 where $a_m(x)$ and $b_m(x)$ are defined in (\ref{1.6}) and (\ref{1.7}). Here $M $ denotes the multiplication operator
 by $|k|$ with domain $D(M) \subset H$. For the third Hamilton equation, one
defines for  all $(q, p)$ in $H^2$ an operator $T(q , p) $ in ${\cal L}( {\cal H}_{sp})$  by,
%----
\be\label{2.8}T (q, p) = \sum _{\lambda = 1}^N \sum _{m=1}^3 ( \beta _m + B_m (x_{\lambda }, q, p) )
 \sigma _m^{[\lambda ]}.\ee
 %-----
The element $i T(q , p)$ belongs to ${\cal G}$. Thus, the Hamiltonian function is written as,
%--
\be\label{2.9} H(q, p, \omega ) = H_{ph} (q , p) +h   <\omega , i T(q, p)  >,\qquad (q, p)\in H^2, \ \ \ \ \omega  \in \Omega.\ee
%----
The third Hamilton equation is then written as, for all $Y\in {\cal G}$, with  (\ref{2.2})
%------
$$ {d \over dt }  <\omega (t), Y> =   {d \over dt }  \phi _Y (\omega (t)) =
\{ H , \phi _Y \} (\omega (t)) = \{ h \phi _{ i T(q(t) , p(t))}, \phi _Y \} (\omega (t)) $$
%----
$$ = \phi _{ [i T(q(t) , p(t)), Y] }  (\omega (t)). $$
%----
In other words, one should have, for all $Y$ in ${\cal G} $,
 %---
\be\label{2.10} {d \over dt }  <\omega (t), Y> =  \sum _{\lambda = 1}^N \sum _{m=1}^3
 ( \beta _m + B_m (x_{\lambda }, q(t), p(t)) )  < \omega (t)) , \Big [ i \sigma _m^{[\lambda ]}, Y
 \Big ] >.\ee
 %-------

Set, for all $(q, p)$ in $H^2$ and for all $t\in \R$,
%----
\be\label{2.11} \chi_t (q , p) = (q_t , p_t),\qquad
\left \{ \begin{matrix}
q_t (k) &=& \cos (t |k|) q(k) +  \sin (t |k|) p(k)\\
 p_t (k) &=& -  \sin (t |k| ) q(k) +  \cos (t |k|) p(k)  \end{matrix} \right. .\ee
%----
We say that a function $t\mapsto ( q(t), p(t), \omega (t)) = ( X(t) ,\omega (t))  $, continuous from $\R$ into
$H \times H \times {\cal G}^{\star}$, is a weak solution to the Hamiltonian system if,
%-----
$$ X(t) = \chi_t (X (0)) + \int _0^t \chi_{t-s} (Q(s), P(s)) ds, $$
%----
$$  Q(s) = h  \sum _{\lambda = 1}^N \sum _{m=1}^3
 b_m (x_{\lambda } )   \ < \omega (s)) , i \sigma _m^{[\lambda ]}>,\qquad
 P(s) = - h \sum _{\lambda = 1}^N \sum _{m=1}^3
 a_m (x_{\lambda }  )   \ < \omega (s)) , i \sigma _m^{[\lambda ]}>  $$
 %----
and if, for all $Y$ in ${\cal G} $,
%----
$$<\omega (t), Y> = <\omega (0), Y>  +   \sum _{\lambda = 1}^N \sum _{m=1}^3 \int_0^t
 ( \beta _m + B_m (x_{\lambda }, q(s), p(s)) )  < \omega (s)) , \Big [ i \sigma _m^{[\lambda ]}, Y
 \Big ] >  ds. $$
 %---

We call a strong solution, a map $t\mapsto ( q(t), p(t), \omega (t)) $
continuous from $\R$ to $D(M) \times D(M) \times  {\cal G}^{\star}$,  $C^ 1$ from $\R$ to
$ H \times H \times  {\cal G}^{\star}$, satisfying (\ref{2.6})(\ref{2.7}) and (\ref{2.10}).

The set ${\cal G}$ is equipped with the scalar product,
%---
\be\label{2.12} X \cdot Y = - {\rm Tr } (X Y),\qquad X \in {\cal G}, \ \ \ \ Y \in {\cal G}\ee
%---
and the dual ${\cal G}^{\star}$ is endowed with the associated Euclidean structure.

\begin{theo}\label{t2.1} For all $(q_0, p_0, \omega _0)$ in $H \times H \times {\cal G}^{\star}$,
there exists a unique weak solution $(q(t), p(t), \omega (t))$ on $\R$ to the system, satisfying $q(0) = q_0$, $p(0)= p_0$ and $\omega (0) = \omega _0$.  One has,
%---
\be\label{2.13} |\omega (t) | =  |\omega (0) |,\qquad |X(t)| \leq |X(0)| + C |t| |\omega (0)|.\ee
%-----
If $q_0$ and $p_0$ belong to $D(M)$ then
this solution is a strong solution.
\end{theo}

 {\it Proof.} Since the $\chi_t$ is a continuous group according to \cite{S} (also see \cite{P}) then
 for any $(q_0, p_0, \omega _0)$ in $H \times H \times {\cal G}^{\star}$ there is a local in time and unique  weak solution on some interval $[0, T )$. Let us show that  (\ref{2.13}) holds
  on $[0, T )$.
For each $Y$ in ${\cal G}$, we denote by ${\rm ad} ^{\star} Y $ the adjoint mapping of  $ X \mapsto {\rm ad}Y ( X)
= [ Y , X ]$.  Equality (\ref{2.9}) is written as,
%---
$$ \omega '(t) = {\rm ad} ^{\star}( i^{-1} T(q (t), p(t)) ) \omega (t).$$
%---
Since the mapping ${\rm ad}( i^{-1} T(q , p)) $ is antiselfadjoint, this is also true for
 ${\rm ad} ^{\star}( i^{-1} T(q , p)) $. Besides, the mapping  $\chi_t$
 is unitary. Inequality (\ref{2.13}) then follows. According to Pazy \cite{P} (Corollary 2.3) or \cite{B-dP-F},
 one concludes that $T = + \infty $. If $q_0$ and $p_0$ belongs to $D(M)$ then
the solution is  strong by \cite{B-dP-F}.

\hfill$\Box$

\begin{prop}\label{p2.2}
If $( q_1(t), p_1(t), \omega _1(t))$ and $( q_2(t), p_2(t), \omega _2(t))$  are two weak solutions then one has,
%----
\be\label{2.14} | X_1(t) - X_2 (t) | + |\omega _1(t) - \omega_2 (t) | \leq A \Big [ | X_1(0) - X_2 (0) |
 + |\omega _1(0) - \omega_2 (0) | \Big ] e^{B |t|},\ee
 where $A$ is a constant and $B$ a function of $|\omega _1(0) | + |  \omega_2 (0) |$.
\end{prop}

{\it Proof.} One first observes that,
%----
$$ | X_1(t) - X_2 (t) | \leq | X_1(0) - X_2 (0) | + A \int_0^t |\omega _1(s) - \omega_2 (s) | ds $$
%---
$$  |\omega _1(t) - \omega_2 (t) |  \leq |\omega _1(0) - \omega_2 (0) | + C \int_0^t
 | X_1(s) - X_2 (s) |ds, $$
 %---
where $C$ is a function of $|\omega _1(0) | + |  \omega_2 (0) |$. Inequality (\ref{2.14}) then follows by Gronwall Lemma.

\hfill$\Box$

\begin{prop}\label{p 2.3}  Let
$(q(t), p(t), \omega (t))$  be a weak solution to the system.
Then $\omega (t)$ lies in the orbit of $\omega (0)$ of the coadjoint representation.
 \end{prop}

 {\it Proof.} There exists a function $g(t)$ taking values in
 $ SU( {\cal H}_{sp})$ satisfying,
 %----
 $$ g'(t) = i T ( q(t), p(t))  g(t),\qquad
 g(0) = I, $$
 %----------
 where $T(q, p)$ is defined in (\ref{2.8}). One sees that, for all $Y$ in ${\cal G} $,
 %---
 $$ < \omega (t), Y> = < \omega (0) , g(t )^{\star } Y g(t) >. $$
 %---
 Consequently, $\omega (t) $ belongs to the  coadjoint orbit of $\omega (0)$.

 \hfill$\Box$

Since the solution $(q(t), p(t), \omega (t))$, a priori belonging to $H^2 \times {\cal G}^{\star}$, remains in our  symplectic manifold then system (\ref{2.6}) (\ref{2.7})(\ref{2.10}) is now called
{\it Hamiltonian system}. We associate  a differential   equation  in  ${\cal H}_{sp}$ with any weak solution $( q(t), p(t), \omega (t))$ to the Hamiltonian system,
 %----
\be\label{2.15} a'(t) = - i T( q(t), p(t) )   a(t)  +   i < T( q(t), p(t) )  a(t) , a(t) > { a(t) \over |a(t)|^2 }.\ee
  %----
  When $a(0)$ belongs to the unit sphere of ${\cal H} _{sp}$, one sees that
  $a(t)$ remains in that sphere and if  $\pi _{a(0)} = \omega (0)$, one has
    $\pi _{a(t)} = \omega (t)$, for all $t$.

\section{Approximate evolution of coherent states.}\label{s3}

We shall derive Theorem \ref{t1.1}. Note that in view of applications,  the norm of $u_h(t) - v_h(t)$ has to be estimated, where $u_h(t)$ and $v_h (t)$ are defined in (\ref{1.18}) and (\ref{1.20}),
not only in the ${\cal H}_{ph} \otimes {\cal H}_{sp} $ norm, but also in a Sobolev spaces norm related to the number operator  $N$.

We denote by  $W_m$ the domain of the
operator $N^{m/2}$, endowed with the norm,
%-----
\be\label{3.1} ||u||_{W_m} ^2=   <(I +2h  N )^m u , u >.\ee
%---
We shall prove the result below which is stronger than Theorem \ref{t1.1}.

\begin{theo}\label{t3.1}  Let $H(h)$ be the Hamiltonian defined in (\ref{1.14})(\ref{1.15}).
Set $X_0 = (q_0 , p_0)$ in $H^2$. Let $a_0$ belongs to ${\cal H}_{sp}$ with norm $1$. Let $(q(t) , p(t) , \omega (t) )
= (X(t) , \omega (t) ) $ be the weak solution to the Hamiltonian system (\ref{2.6})(\ref{2.7})(\ref{2.10}) with   initial data
$(q_0 , p_0, \pi_{a_0})$.
Let  $a(t)$ be the solution to the differential system (\ref{2.15})  with  initial data $a_0$.
Let $u_h(t)$ and $v_h (t)$ be the elements defined in (\ref{1.18}) and (\ref{1.20}).

Then
%---
\be\label{3.2} \Vert u_h(t) - v_h(t) \Vert_{W_2\otimes {\cal H}_{sp}}  \leq C(t) h^{1/2}.\ee
%---
The above constant $C(t)$ is bounded on all compact sets of $\R$.
\end{theo}

The proof uses the five following Lemmas.

\begin{lemm}\label{l3.2} Let $t \mapsto X(t) = (q(t), p(t))$ be a $C^1$ map from $\R$ to $H^2$. Then, the function $t \mapsto \Psi_{q(t), p(t), h}$ is
$C^1$ from $\R$ into ${\cal H}_{ph}$ and  using notation (\ref{1.8}),
%----
\be\label{3.3} h {\partial \over \partial t} \Psi_{q(t), p(t), h} =  Op_h^{weyl} (F_{ q'(t), p'(t)}) \Psi_{q(t), p(t),  h} +
 \gamma (t) \Psi_{q(t), p(t), h},\ee
%------
with
%---
\be\label{3.4}\gamma (t) = - \Big ( q(t) \cdot q'(t) + p(t) \cdot p'(t) \Big ) + {i \over 2}
\Big (  q(t) \cdot  p'(t) - p(t) \cdot q'(t)      \Big ).\ee
%---
\end{lemm}

{\it Proof.}
The  $C^1$ property comes from \cite{D-G} (Proposition 2.4 (iii)). See also Theorem
X.41 (d) of \cite{RE-SI}.
We identify ${\cal H}_{ph}$  with $L^2(B, \mu_{B , h/2})$.
In \cite{A-J-N} or in \cite{A-L-N-1} (Section 2.2), a function $\ell _a $ defined almost everywhere on $B$ is associated with any
element $a$ in the complexified space of  $H$, satisfying  for all $(q, p)$ in $H^2$,
%-----
\be\label{3.5}  \Psi_{q, p , h} (u) = e^{{1\over h} \ell _{ (q+ip)} (u)  -{1\over 2h}|q|^2 -
{i\over 2h} q\cdot p},\qquad   {\rm a.e.\ } u\in B\ee
%---
and such that,  for all $A$ and $B$ in $H$, one has  with  the notation (\ref{1.8}), (\cite{A-J-N}, formula (137)),
%---
\be\label{3.6} (Op_h ^{weyl} F_{A, B} ) = \ell _{  A+i B } (u) + { h\over i } B \cdot {\partial \over
\partial u},\ee
%-----
where  $\ell _{  A+i B } (u)$ stands for the multiplication by the function $\ell _{  A+i B } (u)$.
The derivation a priori has a meaning only on cylindrical functions but  one may extend the equality and shows that,
%---
\be\label{3.7} (Op_h ^{weyl} F_{A, B} ) \Psi_{q, p , h} (u) = \Big [ \ell _{  A+i B } (u) + {1\over i} B \cdot (q+ip)
\Big ] \Psi_{q, p , h} (u).\ee
%---
The Lemma then follows.

\hfill $\Box$

Let us recall (\cite{RE-SI}) that, for all $X$ in $H^2$, there exists a unitary transform $V(X)$
(depending on $h$) satisfying,
%---
\be\label{3.8}\Psi_{X,h}=V(X)\Psi_{0,h}\ee
  %-----
and such that, for all continuous linear form  $L$ on $H^2$,
%----
\be\label{3.9} V(X)^{\star} \Big ( Op_h^{weyl} (L) - L(X) \Big ) V(X)= Op_h^{weyl}(L)\ee
%-----
and verifying, for all $X = (q , p)$ in $D(M)^2$,
%----
\be\label{3.10} V(X)^{\star}  H_{ph} V(X) =  H_{ph} + Op_h^{weyl} ( F_{MX}) + H_{ph } (X).\ee
%--------

\begin{lemm}\label{l3.3} Let $M: H \rightarrow H$ be the multiplication operator by $\omega (k)=|k|$
with $D(M) \subset H$.
Then one has,
%----
\be\label{3.11}  H_{ph}  \Psi_{X, h}  =  Op_h^{weyl} ( F_{MX})  \Psi_{X, h}  - H_{ph } (X)\Psi_{X, h}.\ee
%--------
In addition, for all $(q, p)$ in $D(M)^2$,
%----
\be\label{3.12} [ H_{ph} , Op_h^{weyl} (F_{ q, p} ) ] = i h Op_h^{weyl}
( F_{ - Mp, Mq}).\ee
%----------
\end{lemm}

{\it Proof.}  Equality (\ref{3.11}) follows from (\ref{3.8})(\ref{3.10})
 and (\ref{3.9})  applied to $L = F_{MX} $ together with the fact that $ H_{ph} \Psi_{0, h} = 0$
and  $F_{MX} (X) = 2 H_{ph}(X)$. Equality (\ref{3.12}) is standard.

\hfill $\Box$

\begin{lemm}\label{l3.4} Let $F$ be a continuous linear form  on $H^2$ satisfying $F( -p, q) = i F(q, p)$ for all $(q, p)$ in $H^2$.  Then, one has for all $X$ in $H^2$,
%----
\be\label{3.13} Op_h^{weyl} (F) \Psi _{X h} = F(X) \Psi _{X h}.\ee
%-----
\end{lemm}

{\it Proof.} Let $V(X)$ be the unitary operator  satisfying (\ref{3.8}) and (\ref{3.9}).
In particular,
%------
\be\label{3.14}  Op_h^{weyl} (F) \Psi_{X,h} = F(X) \Psi_{X,h} + V(X)Op_h^{weyl}(F) \Psi_{0,h}.\ee
%----
Under the hypothesis of the Lemma, one has  $Op_h^{weyl}(F) \Psi_{0,h} = 0$ since $Op_h^{weyl}(F)$
is an annihilation operator.

\begin{lemm}\label{l3.5}  For any $(a , b) $ and each $X $ in $H^2$, one has,
%---
\be\label{3.15} \left  \Vert \Big ( Op_h^{weyl} (F_{a , b})   - F_{a , b} (X) \Big  )
 \Psi _{X, h} \right \Vert_{W_2}  \leq C (1+ |X|)^2  h^{1/2}.\ee
%-----
\end{lemm}

{\it Proof.} One uses the operator $V(X)$ satisfying  (\ref{3.8}) and (\ref{3.9}).
It satisfies,
%------
\be\label{3.16}  Op_h^{weyl} (F) \Psi_{X,h} = F(X) \Psi_{X,h} + V(X)Op_h^{weyl}(F) \Psi_{0,h}.\ee
%----
One has, for all $u$ in $W_2$
%---
$$ \Vert u \Vert _{X_2} \leq \Vert u \Vert  + 2h \Vert N u \Vert. $$
%----
It is sufficient to show that,
%----
\be\label{3.17}2h \Vert N  V(X)Op_h^{weyl}(F) \Psi_{0,h} \Vert \leq C (1+ |X|)^2  h^{1/2}.\ee
%------
One has,
%----
$$ 2h N=\sum_k  Op_h^{weyl} (A_{2,k})Op_h^{weyl} (A_{1,k})$$
%---
with
%---
$$ A_{1,k}(q,p)=<q+ip,e_k>,\qquad A_{2,k}(q,p)=<e_k,q+ip>.$$
%------
One applies (\ref{3.9}) with  $L=A_{1,k}$ and then with
 $L=A_{2,k}$ to estimate the left hand side of (\ref{3.17}).
 One gets,
 %----
$$2h N V(X)Op_h^{weyl}(F)=\sum_k
 A_{1,k}(X) A_{2,k}(X) V(X)Op_h^{weyl}(F)
+A_{1,k}(X)V(X) Op_h^{weyl}(A_{2,k})Op_h^{weyl} (F) $$
%----------------------------------------------------
$$
+A_{2,k}(X)V(X)Op_h^{weyl}(A_{1,k})Op_h^{weyl}(F)
+V(X)Op_h^{weyl} (A_{2,k})Op_h^{weyl} (A_{1,k})Op_h(F).$$
%------------------------
Since $Op_h^{weyl} (A_{1,k})$ is an annihilation operator, one has
%---
$$ Op_h^{weyl} (A_{1,k})Op_h^{weyl}(F) \Psi _{0, h} = \Big [ Op_h^{weyl} (A_{1,k}) , Op_h^{weyl}(F)\Big ]  \Psi _{0, h}
 = {h\over i} \{ A_{1, k} , F \} \Psi_{0, h}. $$
 %-----------
One sees  $\sum_k A_{1,k}(X) A_{2,k}(X) =|X|^2/(2h)$. Set
%---
$$ G_X (q , p) = \sum _k \Big ( A_{1,k}(X) A_{2,k}(q, p ) + A_{2,k}(X) A_{1,k}(q, p ) \Big ), $$
%-----
$$ H_X (q , p) = {1\over i} \sum _k  \{ A_{1,k} , F \}  A_{2,k}(q, p ).$$
%---
One has,
%----
$$2h N V(X)Op_h^{weyl}(F) \Psi _{0, h} = |X|^2 V(X)Op_h^{weyl} (F)\Psi _{0, h}   +
  V(X) Op_h^{weyl}( G_X)  Op_h^{weyl} (F)\Psi _{0, h} $$
  %---
  $$+ h V(X) Op_h^{weyl}( H_X) \Psi _{0, h}. $$
%-----
One classically has,
%----
$$ \Vert Op_h^{weyl} (F)\Psi _{0, h} \Vert \leq \sqrt{2}|F| h^{1/2},  $$
%----
$$ \Vert Op_h^{weyl}( G_X) Op_h^{weyl} (F)\Psi _{0, h} \Vert \leq 4h  |G_X| |F|. $$
%---
Consequently,
%---
$$ 2h \Vert  N V(X)Op_h^{weyl}(F) \Psi _{0, h} \Vert \leq  |X|^2 \sqrt{2}|F| h^{1/2}
+ 4h  |G_X| |F|  + \sqrt{2} h^{3/2} |H_X|. $$
%----
One then deduces (\ref{3.17}) and then the  Lemma follows.

\hfill $\Box$

\begin{lemm}\label{l3.6}  There exists $C>0$ such that, for all $f$ in $W_1$ and for any $(a , b)$ in $H^2$,
%---
\be\label{3.18} \Vert Op_h^{weyl} (F_{a , b }) f \Vert \leq C ( |a| + |b|) \Vert f \Vert _{W_1}.\ee
%------
One has for all $f$ in $W_m\otimes {\cal H}_{sp} $ and for any $t\in \R$,
 %----
\be\label{3.19} \Vert e^{ - {i\over h} t H (h) } f \Vert_{W_m\otimes {\cal H}_{sp}}
 \leq C(t)  \Vert  f \Vert_{W_m\otimes {\cal H}_{sp}}.\ee
 %---
\end{lemm}

{\it Proof.}  The first inequality is standard. For the second one, one
uses the following factorization of \cite{A-L-N-2} (Theorem 1.1),
 %----
 $$ e^{ - {i\over h} t H (h) } = e^{-i{t \over h}   H_0} U_h^{red} (t), $$
 %----
 where $H_0 = H_{ph} \otimes I$ and $H_{ph}$ is the free photons Hamiltonian operator.
 Since $H_{ph}$ commutes with  the number operator, one gets for all $f$ in $W_m$,
 %----
 $$ \Vert e^{ - {i\over h} t H _0 } f \Vert_{W_m\otimes {\cal H}_{sp}}
 \leq C(t)  \Vert  f \Vert_{W_m\otimes {\cal H}_{sp}}.$$
 %---
 One now uses the Weyl calculus developed in \cite{A-J-N} and \cite{A-L-N-1}. It is shown
 in \cite{A-L-N-2} that $ U_h^{red} (t)$ is associated by this Weyl calculus with a function
 $H^2 \ni (q , p) \mapsto U(t, q , p, h)) $ taking matrices values (in
${\cal L} ({\cal H }_{sp} ))$. As a function of $(q , p)$, it belongs to in the space
$S_{\infty}^{mat} ({\cal B}, |t| \varepsilon (t))$ introduced in Definition 3.1
of \cite{A-L-N-2} associated with a
Hilbertian basis ${\cal B}$ of $H$ constructed in Section 2.2 of \cite{A-L-N-2}
and with a sequence $\varepsilon (t) $ defined in  Proposition 5.1 of \cite{A-L-N-2}.
In addition, this function is bounded in this space (in the sense of Definition 3.1) as  $h$
is running over $(0, 1)$ and $t$ in  compact sets of $\R$. According to Proposition 2.8 of
\cite{A-L-N-1}, these properties imply that,
%---
$$ \Vert U_h^{red} (t)  f  \Vert_{W_m \otimes {\cal H}_{sp} }  \leq C(t) \Vert   f
  \Vert_{W_m \otimes {\cal H}_{sp}},$$
  %--
  where the constant $C(t)$ is bounded on all compact sets of $\R$. Proposition
  2.8 actually proves the analog of this inequality for the space $W_1$ but the proof for
  $W_m$ is identical.  Inequality  (\ref{3.19}) then follows.

\hfill $\Box$

 {\it End of the proof of Theorem \ref{t3.1}.} Let $(q(t), p(t) , \omega (t))$ be a stong solution to
the Hamiltonian system (\ref{2.6})(\ref{2.7})(\ref{2.10}),
with   $(q(0), p(0), \omega (0))$ in $D(M) \times D(M) \times \Omega $. Set $z_h(t)$ the function
taking values in ${\cal H}_{ph}$ defined by,
%---
\be\label{3.20} z_h(t) = \Psi_{ q(t), p(t), h }.\ee
%--------
From  Lemma \ref{l3.2}, one may  write,
%----
$$ h  {\partial z_h (t)\over \partial t} =  Op_h^{weyl} (F_{ q'(t), p'(t)}) z_h (t) +
 \gamma (t) z_h (t) $$
%------
where $\gamma (t)$ is defined by (\ref{3.4}). From  Lemma \ref{l3.3} one sees,
%----
$$  H_{ph}  z_h (t) =  Op_h^{weyl} ( F_{Mq(t), Mp(t)})   z_h (t)  -
  H_{ph } (q(t), p(t))  z_h (t).  $$
%--------
From (\ref{2.6}) and (\ref{2.7}),
%-----
$$  F_{-p'(t), q'(t)} =  F_{Mq(t), Mp(t)} + h  \sum _{\lambda = 1}^N \sum _{m=1}^3
 B_m (x_{\lambda }, \cdot )   <\omega (t) , i  \sigma _m^{[\lambda ]}  >.
$$
%--------
The function $G = F_{q'(t), p'(t)}  + i F_{-p'(t), q'(t)}  $
satisfies the hypothesis of Lemma \ref{l3.4}.  Therefore,
%----
$$ Op_h^{weyl} (G) z_h (t) = G(q(t), p(t)) z_h (t). $$
%-----
In consequence,
%----
$$  h  {\partial z_h (t)\over \partial t} + i H_{ph} z_h (t) = -i h  \sum _{\lambda = 1}^N \sum _{m=1}^3
   <\omega (t) , i  \sigma _m^{[\lambda ]}  >  B_m (x_{\lambda }, \cdot )  z_h (t)
   + i \lambda (t)  z_h (t) $$
   %---
   with, according to the notation (\ref{1.8}),
   %---
   $$ i \lambda (t)= \gamma (t) + F_{q'(t), p'(t)} (q(t), p(t) )   + i F_{-p'(t), q'(t)}
    (q(t), p(t) )  - i H_{ph} ( q(t), p(t)).  $$
%-------------------
From (\ref{3.4}), one obtains,
%----
$$ \lambda (t) = {1\over 2} \Big (  p(t) \cdot q'(t) -  q(t) \cdot p'(t)\Big )
- H_{ph} ( q(t), p(t)).  $$
%-------------------
Set $w_h (t) = e^{ {i\over h} \varphi (t)} z_h(t)$ where $\varphi (t)$ is defined in (\ref{1.19})
and $z_h(t)$ is defined in (\ref{3.20}).
One gets from (\ref{1.19}),
%----
$$ \varphi '(t) = {1\over 2} \Big (  p(t) \cdot q'(t) -  q(t) \cdot p'(t)\Big )
- H ( q(t), p(t), \omega (t) )  $$
%----
where $H ( q, p, \omega  )$ is the total Hamiltonian function defined in (\ref{2.5}).
Consequently,
%---
\be\label{3.21}h  {\partial w_h (t)\over \partial t} + i H_{ph} w_h (t) = -i h  \sum _{\lambda = 1}^N \sum _{m=1}^3
   <\omega (t) , i  \sigma _m^{[\lambda ]}  >  B_m (x_{\lambda }, \cdot )  w_h (t) + i \mu (t) w_h (t)
\ee
   %---
\be\label{3.22}  \mu (t) =  p(t) \cdot q'(t) -  q(t) \cdot p'(t) - H_{ph} ( q(t), p(t)) -
   H ( q(t), p(t), \omega (t) ).\ee
   %-----
  From (\ref{1.4})(\ref{2.5})(\ref{2.6}) and (\ref{2.7}),
   %---
\be\label{3.23}  \mu (t) =  -h \sum _{\lambda = 1}^N \sum _{m=1}^3 \beta _m
   <\omega (t) , i \sigma _m^{[\lambda ]}  >.\ee
   %-------
Let $a(t)$ be the solution to (\ref{2.15}) with  initial data $a_0$. Let
%----
$$v_h (t) = w_h (t) \otimes a(t)= e^{ {i\over h} \varphi (t)}  \Psi _{q(t) , p(t), h}
\otimes a(t).    $$
%---
Since $\pi _{a_0}  = \omega (0)$, and then $\pi _{a(t)}  = \omega (t)$, one sees from (\ref{1.16}),
%---
\be\label{3.24} <\omega (t) , i \sigma _m^{[\lambda ]}  >  =  < \sigma _m^{[\lambda ]} a(t), a(t) >.\ee
%----
Therefore, equality (\ref{2.15}) together with (\ref{3.21})  give,
%----
$$    h{d \over dt } v_h(t)+ i   H (h) v_h (t)= h S(t, h) $$
%----
with
%---
$$ S(t, h) = - i \sum _{\lambda = 1}^N \sum _{m=1}^3
 \Big (   B_m (x_{\lambda }) - B_m (x_{\lambda }, q(t), p(t) )\Big )
  \Big [   <  \sigma _m^{[\lambda ]} a(t) , a(t) >
   - I \otimes  \sigma _m^{[\lambda ]} \Big ] v_h (t).     $$
%---
Applying Lemma \ref{l3.5} with the function $F(q, p) = B_m (a_{\lambda }, q, p ) $,
one concludes that,
%--
$$ \Vert S(t, h) \Vert_{W_2 \otimes {\cal H}_{sp} } \leq C(t) h^{1/2}, $$
%-----
where the constant $C(t)$ is bounded on all compact sets of $\R$.
One clearly has,
%---
$$ {d \over dt }e^{ - {it\over h} H (h) } \Big ( \Psi_{X_0 , h} \otimes a_0 \Big )
= - {i\over h} H (h) e^{ - {it\over h} H (h) } \Big ( \Psi_{X_0 , h} \otimes a_0 \Big ).$$
%---
Thus,
%---
$$ v_h(t) - e^{ - {it\over h} H (h) } \Big ( \Psi_{X_0 , h} \otimes a_0 \Big ) =
\int _0 ^t e^{ - {i\over h} (t-s) H (h) } S( s, h) ds   $$
%----
From  Lemma \ref{l3.6}, one deduces (\ref{3.2}) if $(q_0, p_0)$ lies in $D(M)^2$.
Now suppose  that $(q_0, p_0)$ only belongs to $H^2$. There exists a sequence $(q_j , p_j)$ in $D(M)^2$  converging
to $(q_0, p_0)$ in $H^2$. Set $(q_j(t), p_j (t) , \omega _j (t))$ the corresponding strong solution  and $a_j (t)$ the solution to (\ref{2.15}). Inequality (\ref{3.2}) is
proved for this solution. In view of
Proposition \ref{p2.2}, the sequence $(q_j(t) , p_j(t))$ tends to $(q(t) , p(t))$ in
$H^2$. In addition, $a_j(t) $ tends to $a(t)$ in ${\cal H}_{sp}$.
From Proposition 2.4 (iii) in \cite{D-G}, $\Psi _{ q_j(t), p_j (t), h }$ tends to
 $\Psi _{ q(t), p (t), h }$ in $W_2$. According to Lemma \ref{l3.6},
 %---
 $$ \lim _{j \rightarrow \infty } \Big \Vert e^{ - {it\over h} H (h) }
 \Big ( \Psi_{q_j, p_j , h} \otimes a_0 \Big ) -
  e^{ - {it\over h} H (h) }
 \Big ( \Psi_{q, p , h} \otimes a_0 \Big ) \Big \Vert _{W_2  \otimes {\cal H}_{sp} } = 0.$$
 %----
 Denoting $\varphi_j (t)$ the action integral analogous to (\ref{1.19}) and
 corresponding to the initial data $(q_j , p_j, \pi _{a_0})$, one also has,
 %---
 $$ \lim _{j \rightarrow \infty } \Big \Vert e^{ {i\over h} \varphi_j (t)}  \Psi _{q_j(t) , p_j(t), h}
\otimes a(t)
 - e^{ {i\over h} \varphi (t)}  \Psi _{q(t) , p(t), h}
\otimes a(t) \Big \Vert _{W_2  \otimes {\cal H}_{sp}  } = 0.$$
 %----
Inequality (\ref{3.2}) then follows for the weak solution.

\hfill $\Box$

\section{Applications. Photons number, fields and spins.}\label{s4}

Our first aim here is to study the average number of photons at time $t$ knowing that the initial state of the system
at $t=0$ is $\Psi_{X_0 , h} \otimes a_0$.

\begin{theo}\label{t4.1} Let $u_h(t)$ be defined in (\ref{1.18}) where $ X_0 = (q_0, p_0) \in H^3$
and let $a_0$ belonging to the unit sphere of ${\cal H}_{sp}$.
Set $(q(t), p(t), \omega (t))$ the solution to the Hamiltonian system with
initial data $(q_0, p_0, \pi _{a_0})$.
Then one has,
%--
$$     \Big | < 2h ( N \otimes I) \  u_h (t) , u_h (t)  > -
  ( q(t)|^2 + |p(t)|^2 ) \Big | \leq C(t)  h^{1/2}, $$
 %----
 where the constant $C(t)$ is bounded on all compact sets of $\R$.

\end{theo}

 {\it Proof.}  Let $u_h (t)$ and   $v_h(t)$ be
 defined in (\ref{1.18}) and (\ref{1.20}). One has,
 %---
 $$ \left | <( N\otimes I ) \  u_h (t) , u_h (t)  > - <( N \otimes I) \  v_h (t) , v_h (t)  > \right | \leq
 2 \Vert ( N\otimes I )( u_h (t) - v_h (t)) \Vert. $$
 %----
 Since the Wick symbol of $N$ is $X = (q, p) \mapsto  (|q|^2 + |p|^2 )/2h$, one sees,
 %---
 $$ < ( N \otimes I) \  v_h (t) , v_h (t)  >  = < N \Psi _{q(t), p(t), h} , \Psi _{q(t), p(t), h} >
 = { |q(t) |^2 + p(t) |^2 \over 2h}.   $$
 %---
Besides, from (\ref{3.1})
 %---
 $$ \Vert ( N\otimes I ) ( u_h (t) - v_h (t) ) \Vert \leq {1 \over 2h } \Vert  u_h (t) - v_h (t) \Vert
 _{W_2 \otimes {\cal H}_{sp} }. $$
 %----
 The Proposition is then completed according to Theorem \ref{t3.1}.

 \hfill $\Box$

Our second aim  is to study the average value of the components of the electric and magnetic fields at each   point $x$ of $\R^3$ and at time $t$,
 knowing  that the initial state of the system at $t=0$ is $\Psi_{X_0 , h} \otimes a_0$.

\begin{theo}\label{t4.2}  Let $u_h (t)$ be
 defined in (\ref{1.18})
where $ X_0 = (q_0, p_0) \in H^3$ and take $a_0$ in the unit sphere of ${\cal H}_{sp}$.
Set $(q(t), p(t), \omega (t))$ the solution to the Hamiltonian system with
initial data $(q_0, p_0, \pi _{a_0})$.  Then  one has,
%--
$$ \Big |  < ( B_m (x)  \otimes I) \  u_h (t) , u_h (t)  > -
  B_m ( x,q(t), p(t)) \Big | \leq C(t) h^{1/2} $$
 %----
where the constant $C(t)$ is bounded on all compact  sets $\R$. This remains also valid  for the electric field.
\end{theo}

{\it Proof.}   Let $u_h (t)$ and   $v_h(t)$ be
 defined in (\ref{1.18}) and (\ref{1.20}). One sees,
 %---
 $$ \left | <( B_m (x)  \otimes I) \  u_h (t) , u_h (t)  > - <( B_m (x)
   \otimes I) \  v_h (t) , v_h (t)  > \right | \leq
2 \Vert (  B_m (x) \otimes I )( u_h (t) - v_h (t)) \Vert.$$
%---
Since the Wick symbol of  $  B_m (x) $ is $B_m (x , p , p)$, one has gets,
 %---
 $$ < ( B_m (x) \otimes I) \  v_h (t) , v_h (t)  >  = < B_m (x) \Psi _{q(t), p(t), h} , \Psi _{q(t), p(t), h} >
 = B_m (x , q(t), p(t).   $$
 %---
Besides, from Lemma \ref{l3.6},
 %---
 $$ \Vert ( B_m (x)\otimes I ) ( u_h (t) - v_h (t) ) \Vert \leq  C \Vert  u_h (t) - v_h (t) \Vert
 _{W_1 \otimes {\cal H}_{sp} }. $$
 %----
 The Proposition is then complete in view of Theorem \ref{t3.1}.

\hfill$\Box$

Our next aim  is to study the average value of the  spins of the particles at time $t$,
 knowing that the initial state of the system at time $t=0$ equals to $\Psi_{X_0 , h} \otimes a_0$.

\begin{theo}\label{t4.3} Let $u_h (t)$ be
 defined in (\ref{1.18})
where $ X_0 = (q_0, p_0) \in H^3$ and let $a_0$ lying in the unit sphere of ${\cal H}_{sp}$.
Let $(q(t), p(t), \omega (t))$ be the solution to Hamiltonian system with
initial data $(q_0, p_0, \pi _{a_0})$.  Let $a(t)$ be the solution to the differential system (\ref{2.15}) with  initial data $a_0$. Then one has, for all $\lambda \leq N$
and for all $m\leq 3$,
%--
$$ \Big |  < ( I \otimes \sigma _m ^{[\lambda]} ) \  u_h (t) , u_h (t)  > -
< \sigma _m ^{[\lambda]} a(t), a(t) > \Big | \leq
    C(t) h^{1/2} $$
 %----
where the above constant $C(t)$ is bounded on all compact sets of $\R$.
\end{theo}

The proof is identical to those of Theorems \ref{t4.1} and \ref{t4.2}.

\section{Approximate evolution laws}\label{s5}

Theorems \ref{t4.1}, \ref{t4.2} and \ref{t4.3} enable to write (semiclassical) approximate laws for
the evolution of the electromagnetic field average, the spin averages and the photons number average.
 We recover classical laws for  fields, laws going back to Bloch \cite{BL} (1946) for
 spins and we state a perhaps new law for the photons number average.

Fix $( q(t), p(t), \omega (t))$ a solution to the  Hamiltonian system and $a(t)$ a
solution to (\ref{2.15}). Set,
\be\label{5.1} B_m ^{appr} (x , t) =  B_m ( x, q(t), p(t)),\qquad
E_m ^{appr} (x , t) =  E_m ( x, q(t), p(t))\ee
%----
and also, using (\ref{3.24}),
%----
\be\label{5.2} S_m ^{[\lambda , appr]} (t) = < \sigma _m ^{[\lambda]} a(t), a(t) >
= < \omega (t), i \sigma _m ^{[\lambda]} >.\ee
%-----
Also set,
 ${\bf B} ^{appr } (x, t) = ( B_1 ^{appr} (x , t) , B_2 ^{appr} (x , t) , B_3 ^{appr} (x , t) )$,
and similarly for ${\bf E} ^{appr } (x, t)$ and $ {\bf S}  ^{[\lambda , appr]} (t)$.

{\it Fields evolution.} First, one remarks that,
%-----
\be\label{5.3} {\rm div} {\bf B} ^{appr } (x, t) = {\rm div} {\bf E} ^{appr } (x, t) = 0.\ee
%----
We shall show that,
%--------
\be\label{5.4} {\partial \over \partial t}  {\bf B} ^{appr}   (x, t)  = - {\rm rot}
 {\bf E} ^{appr}   (x, t)\ee
%------
\be\label{5.5} {\partial \over \partial t}  {\bf E} ^{appr}   (x, t)  =  {\rm rot}
 {\bf B} ^{appr}   (x, t)+ h \sum _{\lambda = 1} ^N {\bf S^{[\lambda , appr]}} (t) \wedge {\rm grad }
 \rho (x-x_{\lambda}),\ee
%------
with,
%---
\be\label{5.6} \rho (x) = (2\pi)^{-3} \int _{\R^3} |\chi (k)|^2
\cos ( k \cdot x   )  dk.\ee
%----
Using definition (\ref{1.11}) of $E_m (x, p, p)$ together with the fact  that $J$ is antisymmetric,
one may write
$E_m (x, p, p) = \alpha _m (x) \cdot q + \beta _m (x) \cdot p$, where $\alpha _m (x)$
and $\beta _m (x)$ are the elements of $H$ defined by,
 %----
\be\label{5.7} \alpha _m (x)   (k) =
 {\chi(|k|)|k|^{1\over 2} \over (2\pi)^{3\over 2}} \sin ( k \cdot x )
{k \wedge (k\wedge e_m) \over |k|^2},\ee
%---
\be\label{5.8} \beta _m (x)   (k) =
 {\chi(|k|)|k|^{1\over 2} \over (2\pi)^{3\over 2}} \cos  ( k \cdot x )
{k\wedge ( k \wedge e_m) \over |k|^2}.\ee
%---
One uses the Hamiltonian system noticing that, for all
$q$ and $p$ in  $D(M)$,
%----
$$ {\bf E} (x , Mp , -Mq) = {\rm rot } {\bf B} (x , q , p),\qquad
{\bf B} (x , Mp , -Mq) = -  {\rm rot } {\bf E} (x , q , p). $$
%----
One then deduces equality (\ref{5.4}) and the fact that,
%------
$$ {\partial \over \partial t}  {\bf E} ^{appr}   (x, t)  =  {\rm rot}
 {\bf B} ^{appr}   (x, t)+ {\bf F} (x , t),$$
 %---
 $$ F_m (x , t ) =  h \sum _{\lambda = 1} ^N \sum _{n = 1} ^3 K _{mn} ( x, x_{\lambda})
 S_n^{[\lambda , appr]} (t),$$
 %----
 with
 %---
 $$ K _{mn} ( x, x_{\lambda})  = \alpha _m(x)   \cdot b_n (x_{\lambda})
 - \beta _m(x)\cdot a_n (x_{\lambda}),   $$
 %----
 where the $a_m(x)$ and $b_m(x)$ are  defined in (\ref{1.6}) and (\ref{1.7}), and $\alpha_m(x)$ and $\beta_m(x)$
 in (\ref{5.7}) and (\ref{5.8}). One then observes that  the matrix $ K_{mn} ( x, x_{\lambda})$ is antisymmetric and that,
for instance,
%---
$$ K _{12} ( x, x_{\lambda})= (2\pi)^{-3} \int _{\R^3} |\chi (k)|^2
\sin ( k \cdot (x - x_{\lambda} ) ) k_3 dk = - ( \partial_3 \rho ) ( x- x_{\lambda })$$
%----
and the two others are obtained by circular  permutation. Thus, one gets  (\ref{5.5}). This approximate evolution law is Maxwell laws,
with  zero charge and a divergence free current modelling the spin.

{\it Evolution of spins.} From  (\ref{5.2}) and (\ref{2.10}),
%---
$$ {d \over dt}  S_m ^{[\lambda , appr]} (t) = {d \over dt }  <\omega (t), i \sigma _m ^{[\lambda]}>
 =  \sum _{\mu = 1}^N \sum _{n=1}^3
 ( \beta _n + B_n (x_{\mu }, q(t), p(t)) )  < \omega (t)) , \Big [ i \sigma _n^{[\mu ]},
 i \sigma _m^{[\lambda ]} \Big ] >. $$
%-------
The commutator vanishes if $\mu \not = \lambda$. Using Pauli matrices commutation formulas, for example $[\sigma _3 , \sigma_1] = 2i \sigma_2$, and  then using again (\ref{5.2}), one obtains,
%---
\be\label{5.9} {d \over dt} {\bf S}  ^{[\lambda , appr]} (t) = 2 ( \beta + {\bf B} ^{appr } (x_{\lambda} , t) )
\wedge {\bf S}  ^{[\lambda , appr]} (t).\ee
%---
These are  Bloch equations \cite{BL} (1946).

 {\it Approximate evolution of the photons number average.}  Set,
   %----
   $$ N ^{appr} (t)= { |q(t) |^2 + p(t) |^2 \over 2h}. $$
   %------
 In view of the Hamiltonian system, this approximate average evolves according to the following law,
   %---
   $$ {d \over dt } N ^{appr} (t)=
 \sum _{\lambda = 1}^N \sum _{m=1}^3
 B_m  (x_{\lambda } , -p(t), q(t)  )
  \ < \omega (t)) , i \sigma _m^{[\lambda ]}>. $$
 %------
  From (\ref{5.2}),
$ < \omega (t)) , i \sigma _m^{[\lambda ]}>  =  S_m ^{[\lambda , appr]} (t)$.
   One may give a more physical meaning to $B_m  (x_{\lambda } , -p(t), q(t)  ) $
   using the helicity operator $J$ from $H^2$ to $H^2$ defined in (\ref{1.10}).
  For all $(q, p)$ in $H^2$, set
 ${\cal F}  (q, p) = (-p, q)$. Since $J$ and  ${\cal F}$ are commuting, one sees  that $J^{-1} {\cal F}$ is
selfadjoint in $H^2$ and that its square equals to the identity. Its eigenvalues are $+1$ and $-1$. Denote by $E_+ $ and $E_-$ the corresponding eigenspaces. The sum of these subspaces  is $H^2$.
 One denotes by $\Pi_+$ and $\Pi_-$ the corresponding  projections. One has for all $Z = (q, p)$,
 %----
 $$ (\Pi_{\pm} Z) (k) = {1\over 2 } \Big [ Z(k) \mp (JZ) (k) \Big ]. $$
 %----
 The subspaces $E_+$ and $E_-$ are respectively corresponding to the right and left circular polarization.
One has  $ {\cal F} \Pi_{\pm} Z = \pm J \Pi_{\pm} Z$. Consequently,
%----
$$ {d\over dt}  N^{appr} (t)   =   \sum _{\lambda = 1}^N \sum _{j=1}^3
( B_j (x_{\lambda} , J \Pi _+(q(t), p(t))  )
-  B_j (x_{\lambda} , J \Pi _- (q(t), p(t) ) ))
  S_j ^{[\lambda , appr]} (t).$$
%-----
Therefore, from  (\ref{1.11}),
%----
\be\label{5.10} {d\over dt}  N^{appr} (t)   =   \sum _{\lambda = 1}^N \sum _{j=1}^3
( E_j (x_{\lambda} , \Pi _- (q(t), p(t) ) ) - E_j (x_{\lambda} ,  \Pi _+(q(t), p(t))  )
)
  S_j ^{[\lambda , appr]} (t).\ee
%-----

\section{Fixed points and quasimodes.}\label{s6}

We shall first give a characterization of the fixed points of the Hamiltonian system.
 For all $u$ in the unit sphere of ${\cal H}_{sp}$, define
 $(q_h (u) , p_h(u))$ as the  point in $H^2$ by,
 %----
\be\label{6.1}   M q _h (u)  = - h \sum _{\lambda = 1}^N \sum _{m=1}^3
 a_m (x_{\lambda }  )   \ <   \sigma _m^{[\lambda ]} u, u>,\ee
 %------
\be\label{6.2} M p _h (u)   = - h  \sum _{\lambda = 1}^N \sum _{m=1}^3
 b_m (x_{\lambda } )   \ <  \sigma _m^{[\lambda ]}u , u >,\ee
 %------
 where $ a_m (x_{\lambda }  )$ and $ b_m (x_{\lambda }  )$ are elements of $H$
 defined in (\ref{1.6}) and (\ref{1.7}).
Since the function $\chi $ in (\ref{1.6}) and (\ref{1.7}) belongs to ${\cal S}(\R)$ then this
 element  $(q_h (u) , p_h(u))$  exists in $H^2$. Even if $\chi$ is not vanishing in a neighborhood of the origin, the function $k\mapsto a_m (x_{\lambda } (k) / |k| $ belongs to $H$.

 \begin{prop}\label{p6.1}
  Let $(q, p)$ be in $H^2$ and $a$ in the unit sphere of
 ${\cal H}_{sp}$. The following properties are equivalent.

 $(i)$ The point  $(q , p, \pi_a)$ is a fixed point  of the Hamiltonian system.

$(ii)$ One has  $q= q_h(a)$, $p= p_h (a)$ and  $a$ is an eigenvector of the operator $T(q , p) $
 defined in (\ref{2.8}).
\end{prop}

 {\it Proof.}  If $(i)$  is satisfied then  equalities (\ref{2.6}) and (\ref{2.7}) show  that
  $q= q_h(a)$ and $p= p_h (a)$.  In view of (\ref{2.10}) one has,
 %---
 $$ < [T(q, p), X ] a, a > = 0,\qquad\forall X \in {\cal G}. $$
 %----
 Set $b$ in ${\cal H}_{sp}$ orthogonal to $a$. There exists $X \in  {\cal G}$
 satisfying $Xa = b$. To check that, it is sufficient to choose an orthonormal basis
 $(e_1 , \dots, e_p)$ ($p= 2^N$) of ${\cal H}_{sp}$ satisfying $a =  e_1$ and $b= \mu e_2$.
 with  $\mu \in \R$ and to choose $X$ such that $Xe_1 = \mu e_2$, $Xe_2 = - \mu e_1$
 and $Xe_j=0$ if $j\geq 3$. One deduces,
 %----
 $$ < [T(q, p), X ] a, a > = 2 {\rm Re} < T(q, p) a , b> = 0.$$
 %---
 Replacing $b$ by $ib$, it is concluded that $< T(q, p) a , b> = 0$ for
 all vectors $b$ orthogonal to $a$. This means that $T(q, p) a $ is colinear
 to $a$ and therefore $a$ is an eigenvector of $T(q, p) $.

\hfill$\Box$

We shall now identify the fixed points of the Hamiltonian system.
The constant magnetic field  $\beta $ is supposed to be non zero. Two unitary vectors
   $b_0$ and $b_1$ of ${\bf C}^2$ are chosen satisfying,
   %----
   %----
\be\label{6.3}\sum _{m=1}^3 \beta _m \sigma_m b_{0 } = - |\beta |
 b_{0},\qquad\sum _{m=1}^3 \beta _m \sigma_m b_{1 } = + |\beta |
 b_{1}.\ee
 %----
For all $E \subset \{ 1 , \dots , N \}$,  $a_E$ denotes the following element,
%---
\be\label{6.4} a_E = a_1 \otimes \cdots \otimes a_N, \qquad  a_j = \left \{ \begin{matrix}  b_1 & {\rm if} &j\in E \\
 b_0 & {\rm if} & j \notin E  \end{matrix}\right ..\ee
 %--
  The eigenvalues of the operator $T(0, 0)$ defined in (\ref{2.8})
are $\mu _p =  (2p -N ) |\beta | $ ($0 \leq p \leq N$). The eigenspaces corresponding to
$ \mu _p$ is spanned by the $a_E$ with  $|E| = p$.

 \begin{prop}\label{p6.2}   The constant magnetic field $\beta $ is assumed to be non vanishing.
 Let $a$ be in  the unit sphere of ${\cal H}_{sp}$ and not belonging to the union of circles  $e^{i \theta } a_E$ ($E \subset \{ 1 , \dots , N \}$, $\theta \in \R$).
 Then, for small enough $h$, there exists no point $(q, p)$ in $H^2$ such that $(q , p , \pi _a)$
 is a stationary point of the Hamiltonian system.
\end{prop}

 {\it Proof.}
 Suppose that  there is a point $(q, p)$ in $H^2$ such that $(q , p , \pi_a)$ is  a stationary point. We necessary have
  $ q = q_h (a)$ and $p = p_h(a)$ (defined in (\ref{6.1}) and (\ref{6.2})).
In view of Proposition \ref{p6.1},
the point $a$ is an eigenvector of $T ( q_h (a) , p_h (a) )$ defined in (\ref{2.8}). Equalities
 (\ref{2.8})(\ref{6.1}) and (\ref{6.2})  show that,
 %---
 $$ \Vert T ( q_h (a) , p_h (a) ) - T (0, 0) \Vert \leq Ch. $$
 %----
Since $a$ is not an eigenvector of $T(0, 0) $ then there is a contradiction for $h$ small enough.

\hfill $\Box$

Conversely, we shall associate fixed points of the Hamiltonian system with elements $a_E$. However, a geometrical hypothesis concerning the position of spin particles is needed.

\begin{prop}\label{p6.3}  The constant magnetic field $\beta $ is assumed to be non vanishing
and all the  particles are supposed to be located at points $x _{\lambda }$
 in the same orthogonal plane to $\beta$.
 Set $a = a_E$ in the unit sphere of ${\cal H}_{sp}$ defined by (\ref{6.3}) (\ref{6.4})
  with  $E \subset \{ 1 , \dots , N \}$.

  Then, for all $h>0$, there exists $(q_h , p_h)$ in
 $H^2$ such that $(q_h , p_h , \pi_a)$ is a fixed point  of the  Hamiltonian system.
Assuming that $\beta = (0, 0, |\beta|)$ and setting $\varepsilon _{\lambda } = 1$ if
$\lambda \in E$ and $\varepsilon _{\lambda } = - 1$ otherwise, the energy level of the fixed point is given by,
%----
%---
\be\label{6.5} H (q_h(a) , p_h (a), \pi_a )= h |\beta | \sum _{\lambda = 1}^N \varepsilon _{\lambda }
 - {h^2 \over 2} \sum _{\lambda , \mu \leq N }    C_{33}^{[\lambda , \mu ]}
 \varepsilon_{\lambda } \varepsilon_{\mu },\ee
 %----------
with
%---
\be\label{6.6} C_{33}^{[\lambda , \mu ]}  = (2 \pi )^{-3}  \int _{\R^3 } |\chi(|k|) |^2
  \cos ( k \cdot (x _{\lambda } - x _{\mu })) \  {k_1^2 + k_2^2  \over |k|^2}
  dk.\ee
 %---
\end{prop}

 {\it Proof.}  One may suppose that $\beta = (0, 0, |\beta|)$.
Set $(q_h (a) , p_h(a))$ the point of  $H^2$ defined by (\ref{6.1}) and (\ref{6.2}).
The operator $T(q_h (a) , p_h(a))$  defined in (\ref{2.8}) is written as,
 %----
 $$ T(q_h (a) , p_h(a)) = |\beta | \sum _{\lambda = 1}^N \sigma _3^{[\lambda ]} -h
 \sum _{\lambda , \mu \leq N }  \sum _{m , n \leq 3  }   C_{mn}^{[\lambda , \mu ]}
 < \sigma _n^{[\mu ]}  a, a> \sigma _m^{[\lambda ]}, $$
 %------
 with
 %---
 $$  C_{mn}^{[\lambda , \mu ]}  = a_m ( x_{\lambda }) \cdot ( M^{-1} a_n ( x_{\mu })) +
 b_m ( x_{\lambda }) \cdot ( M^{-1} b_n ( x_{\mu })), $$
 %----
 where the $a_m(x)$ and $b_m (x)$ are defined in (\ref{1.6}) and (\ref{1.7}).
 Consequently,
%---
 $$  C_{mn}^{[\lambda , \mu ]}  = (2 \pi )^{-3}  \int _{\R^3 } |\chi(|k|) |^2
  \cos ( k \cdot (x _{\lambda } - x _{\mu })) \  {(k\wedge e_m) \cdot (k\wedge e_n) \over |k|^2}
  dk. $$
 %---
  If $a= a_E$ then one has  $< \sigma _n^{[\mu ]}  a, a> = 0$ excepted if $n= 3$.
  If $n= 3$ then one has  $  C_{mn}^{[\lambda , \mu ]}  = 0$ excepted if $m= 3$. Indeed, if $n=3$ and $m\not = 3$,
  it is sufficient to replace $k_3$ by $-k_3$ to see that $  C_{mn}^{[\lambda , \mu ]}  = 0$ since
  $ \cos ( k \cdot (x _{\lambda } - x _{\mu }))$ is independent of $k_3$, since $\chi $
  is radial and since $(k\wedge e_m) \cdot (k\wedge e_3) $ is an odd function of $k_3$
  when $m\not = 3$. Therefore,
  %---
 $$ T(q_h (a) , p_h(a)) = |\beta | \sum _{\lambda = 1}^N \sigma _3^{[\lambda ]} -h
 \sum _{\lambda , \mu \leq N }    C_{33}^{[\lambda , \mu ]}
 < \sigma _3^{[\mu ]}  a, a> \sigma _3^{[\lambda ]}.  $$
 %----
Thus, $a = a_E$ is an eigenvector of $ T(q_h (a) , p_h(a))$ with an eigenvalue being,
 %---
 $$ \lambda _h = |\beta | (2|E| - N)   -h
 \sum _{\lambda , \mu \leq N }    C_{33}^{[\lambda , \mu ]}
 \varepsilon_{\lambda } \varepsilon_{\mu }.   $$
 %----
 Proposition \ref{p6.1} shows that
 $(q_h(a) , p_h (a), \pi_a )$ is a fixed point of the Hamiltonian system.
From (\ref{1.4}) and (\ref{2.5})-(\ref{2.8}),  one has,
 %---
 $$ H (q_h(a) , p_h (a), \pi_a )= {h\over 2 } \Big [ < T(0, 0) a, a > + < T(q_h(a) , p_h (a)) a, a > \Big ]. $$
 %---
Then (\ref{6.5}) is derived.

\hfill$\Box$

In order to give a more physical meaning to  (\ref{6.6}),  $\chi $ has to tend to $1$. The fact
 that $\chi $ is vanishing in a neighborhood of the origin is not used here (it is used in Sections \ref{s3}-\ref{s5}
 and in the below construction of the quasimode).  One may then  take $\chi (|k|) = \chi_{\varepsilon }(|k|) =
  \varphi ( |k|  \varepsilon )$,
 where $\varphi $ is a function $C^{\infty}$ with compact support which equals to $1$ in a neighborhood of $0$.
 Explicitly writing the parameter $\varepsilon$ as an index, one then has,
 %---
  $$ H_{\varepsilon}  (q_h(a) , p_h (a), \pi_a )= h |\beta | \sum _{\lambda = 1}^N \varepsilon _{\lambda }
 - {h^2 \over 2} \sum _{\lambda , \mu \leq N }  F_{\varepsilon } (x_{\lambda } - x_{\mu  } )
 \varepsilon_{\lambda } \varepsilon_{\mu },  $$
 %----------
 $$ F_{\varepsilon } (x ) = (2 \pi )^{-3}  \int _{\R^3 } |\varphi (|k| \varepsilon  ) |^2
  \cos ( k \cdot  x ) \  {k_1^2 + k_2^2  \over |k|^2}
  dk.  $$
 %---
 For all non zero $x$ , one sees that,
 %----
 $$ \lim _{\varepsilon \rightarrow 0 } F_{\varepsilon } (x ) = -
 ( \partial _1^2 +  \partial _2^2 ) { 1 \over 4 \pi |x| }.$$
 %---
 In particular, if $x_3 = 0$,
 %----
 $$ \lim _{ \varepsilon \rightarrow 0 } F_{\varepsilon } (x ) =  - {1 \over 4 \pi |x|^3 }.  $$
 %---
 Set $C_{ \varepsilon} $ the sum of the  terms corresponding to $\lambda = \mu $ in
 (\ref{6.5}). This sum is independent of $E$. One obtains,
 %---
 $$ \lim _{\varepsilon \rightarrow 0 } H_{\varepsilon}  (q_h(a_E) , p_h (a_E), \pi_{a_E} ) -
 C_{ \varepsilon} = h |\beta | \sum _{\lambda = 1}^N \varepsilon _{\lambda }
 + {h^2 \over 8 \pi } \sum _{\lambda \not=  \mu  } { \varepsilon_{\lambda } \varepsilon_{\mu }
\over |x_{\lambda } - x_{\mu  }|^3 }. $$
%---

We shall now associate quasimodes with  points $a_E$ ($E \subset \{ 1 , \dots , N \}$). These quasimodes are of low precision if $E \not = \emptyset$
and of arbitrary high accuracy if  $E= \emptyset$.

\begin{prop}\label{p6.4}
Fix $(q, p)$ in $H^2$ and take $a$ in the unit sphere of
${\cal H}_{sp}$ such that $( q, p, \pi_a)$  is a stationary point of the Hamiltonian system
(\ref{2.6})-(\ref{2.7})-(\ref{2.10}). Then one has,
%---
\be\label{6.7}\Big [  H_{ph} + h \sum _{\lambda = 1}^N \sum _{m=1}^3 ( \beta _m + B_m (x_{\lambda }) )
 < \sigma _m^{[\lambda ]} a , a  > \Big ] \Psi_{q, p, h} = H(q , p, \pi_a )
 \Psi_{q_h, p_h, h},\ee
 %---
 where $ H(q , p, \omega )$ is the Hamiltonian function defined in (\ref{2.5}).

Suppose that $\beta $ is non zero. Set $E \subset \{ 1 , \dots , N \}$.
Then one has, with  a constant $C>0$,
%---
\be\label{6.8}\Vert ( H(h) - (2 |E| - N ) |\beta |h )  ( \Psi_{q, p, h} \otimes a_E) \Vert
\leq C h^{3/2}.\ee
%-----
\end{prop}

 {\it Proof.} From Proposition \ref{p6.1}, if $(q , p, \pi_a)$ is a stationary point
 then $q$ and $p$ are in $D(M)$ and given  by (\ref{6.1}) and (\ref{6.2}).
From The Lemma 3.3, setting $X = (q, p)$,
 %----
$$  H_{ph}  \Psi_{X, h}  =  Op_h^{weyl} ( F_{MX})  \Psi_{X, h}  - H_{ph } (X)\Psi_{X, h}.  $$
%--------
One has, from equalities (\ref{6.1}) and (\ref{6.2}),
 %---
 $$ F_{MX}  =  -  h \sum _{\lambda = 1}^N \sum _{m=1}^3 B_m (x_{\lambda }, \cdot )
  <  \sigma _m^{[\lambda ]} a, a >. $$
  %--------
In view of  (\ref{3.22}) and (\ref{3.23}) (with  $ q' = p' = 0$) and (\ref{3.24}), one sees,
 %--
 $$ - H_{ph } (X) - H(q , p, \pi_a ) = - h \sum _{\lambda = 1}^N \sum _{m=1}^3
 \beta _m  <  \sigma _m^{[\lambda ]} a, a >.  $$
 %----------
 Then (\ref{6.7}) is deduced. For (\ref{6.8}), one may suppose that $\beta = (0, 0, |\beta |)$,
 and that $(q , p) = ( q_h (a_E), p_h (a_E))$.
Then (\ref{6.7}) gives,
 %---
$$\Big [  H_{ph} \otimes I + h \sum _{\lambda = 1}^N  ( \beta _3 + B_3 (x_{\lambda }) )
 \otimes  \sigma _3^{[\lambda ]}   \Big ] \Big ( \Psi_{ q_h (a_E), p_h (a_E) , h} \otimes a_E \Big )
 = H(q_h (a_E), p_h (a_E) ,  \pi_{a_E}  )
 \Big ( \Psi_{q_h (a_E), p_h (a_E) , h} \otimes a_E \Big ).  $$
 %---
 Consequently,
 %----
 $$ ( H(h) - H( q_h (a_E), p_h (a_E) , \pi_{a_E}  ) ) \Big ( \Psi_{q_h (a_E), p_h (a_E) , h} \otimes a_E \Big )  =
h \sum _{\lambda = 1}^N \sum _{m=1}^2  B_m  (x_{\lambda }) \otimes  \sigma _m^{[\lambda ]}
\Big ( \Psi_{q_h (a_E), p_h (a_E) , h} \otimes a_E \Big ). $$
%----
Since $|  B_m  (x_{\lambda }, q_h (a_E), p_h (a_E) )| \leq Ch$ then  Lemma \ref{l3.5} shows that the
norm in the right hand side is $\leq C h^{3/2}$. Besides, $H( q_h (a_E), p_h (a_E) , \pi_{a_E}  ) =
(2 |E| - N ) |\beta |h  +  {\cal 0} (h^2)$. Thus (\ref{6.8}) is derived.

\hfill $\Box$

The following theorem provides additional precisions when the set  $E$ is
empty.

\begin{theo}\label{t6.5}  Suppose that $\beta $ is non zero and that the function $\chi$ in (\ref{1.6}) and (\ref{1.7})
is vanishing in a neighborhood of $0$. With the notations (\ref{6.3}) (\ref{6.4}),
set $a= a_{\emptyset} \in {\cal H}_{sp}$. Then, there exists
 a sequence of elements of ${\cal H}_{ph} \otimes {\cal H}_{sp} $  denoted
 $u_{j}$ ($j\geq 0$) together with a sequence of real numbers $\lambda _j$ ($j\geq 1$) satisfying,
 %----
\be\label{6.9} u_0 = \Psi_{0} \otimes a_{\emptyset},\qquad \lambda _1 = - N |\beta | \ee
 %---
 and such that,
%---
\be\label{6.10} \left \Vert \Big (H(h) -(\lambda _1 h + \cdots + \lambda _{p+1} h^{p+1} ) \Big ) \ \left [ \sum _{j=0}^{2p}
 u_{j} h^{j\over 2} \right ] \right \Vert \leq C  h^{ p + (3/2)},\ee
 %-----
\be\label{6.11} \left \Vert \Big (H(h) -(\lambda _1 h + \cdots + \lambda _{p+1} h^{p+1} ) \Big ) \ \left [ \sum _{j=0}^{2p+1}
 u_{j} h^{j\over 2} \right ] \right \Vert \leq C  h^{ p + 2}.\ee
 %-----
\end{theo}

The above elements $u_{j}$  are independent of $h$ provided that these are considered as element of the Fock space ${\cal F}_s(H_{\bf  C})$ and not of
$L^2 (B, \mu _{B , h/2})$. Indeed, the isomorphism between these two spaces depends on $h$.

{\it Additional details on  Fock space.} In Sections \ref{s2} to \ref{s5}, we have considered that ${\cal H}_{ph}$
is a  $L^2 (B, \mu _{B , h/2})$ space associated with a  Wiener space
related to $H$. From now on,  it is more convenient to suppose that it is the symmetrized Fock space
%----
\be\label{6.12} {\cal F}_s (H_{\bf C}) =  \oplus _{m\geq 0} {\cal F}_m, \ee
%--
where ${\cal F}_0 = {\bf C}$ and ${\cal F}_m$ is completion of the $m-$fold symmetric tensor product   $H_{\bf C} \odot \cdots \odot H_{\bf C} $. One may
then consider an element of ${\cal F}_m$ as a symmetric map  $f$ from
$ (\R^3)^m $ to $ (\R^3) ^{\otimes m}$ satisfying for all $a_2,\dots,a_m$ in $\{ 1, 2, 3 \}$ and for all $k_1,\dots,k_m$ in $\R^3$,
%---
\be\label{6.13} \sum _{j=1}^3 { k_{1 , j}}  f_{ j, a_2 , ... a_m} ( k_1 , \dots , k_m) = 0.\ee
%---
We use here the notation $k_1 = (  k_{1 , 1} ,  k_{1 , 2} ,  k_{1 , 3})$. In addition, the components of this function $f$ should be in $L^2 (\R^{3m})$, which is  defining the norm in ${\cal F}_m$.
Thus, an element of $ {\cal F}_s (H_{\bf C})$ is a sequence $ f = (f_m)_{(m\geq 0 )}$ where $f_m$ is an
element of ${\cal F}_m$ and one has,
%---
\be\label{6.14} \Vert f \Vert ^2 = \sum _{m\geq 0} \Vert f_m \Vert ^2.\ee
%---

For any $\rho >0$ and $m\geq 1$,  ${\cal F}_m (\rho )$ stands for the set of elements
$f$ in  ${\cal F}_m  $ satisfying $f(k_1 , ... , k_m) =0$ if one of the $|k_j|$ is
$\leq \rho $. If $m=0$, it is agreed that ${\cal F}_0 (\rho )= {\cal F}_0 $. It is also agreed
that ${\cal F}_m   =0$ if $m<0$. One sets,
%---
\be\label{6.15} {\cal F} _{even} (\rho) =  {\cal F}_0 \oplus {\cal F}_2 (\rho) \oplus {\cal F}_4 (\rho) \oplus \cdots\ee
%---
\be\label{6.16} {\cal F} _{odd} (\rho) =  {\cal F}_1 (\rho)  \oplus {\cal F}_3 (\rho) \oplus {\cal F}_5 (\rho) \oplus \cdots.
\ee
%---
These elements here are finite sums.

Let us recall that, if $M: H \rightarrow H$ is the multiplication by
$\omega (k) = |k|$ and if $f$ is a rapidly decreasing function in ${\cal F}_m$
then one has,
%---
\be\label{6.17} ( {\rm d}\Gamma (M) f ) ( k_1 , \dots , k_m) = ( |k_1| + \dots + |k_m | ) f  ( k_1 , \dots , k_m).\ee
%----
If $L(q , p) = a \cdot q + b \cdot p$ is a continuous linear form  on $H^2$, we associate with it,
for all $h>0$, an unbounded operator in  $L^2 (B, \mu _{B , h/2})$ denoted $Op_h^{weyl} (L)$.
This operator is defined for smooth cylindrical functions  by (\ref{3.6}). Via
the  canonical (Segal) isomorphism
 ${\cal J}_h : {\cal F}_s (H_{\bf C}) \rightarrow L^2 (B, \mu _{B , h/2})$, this operator becomes,
 %----
\be\label{6.18} {\cal J}_h ^{-1} Op_h^{weyl} (L)  {\cal J}_h = \sqrt {h}  \Phi_S (a+i b),\ee
 %----
 where $\Phi_S (a+i b) :{\cal F}_m \rightarrow {\cal F}_{m+1} \oplus {\cal F}_{m-1} $ is a continuous linear operator,  called Segal field and defined in  \cite{RE-SI}.
 We denote by $\Psi_0$ a unitary element of ${\cal F}_0$. One may suppose that
 ${\cal J}_h  \Psi_0 = \Psi_{0, h}$ where $\Psi_{0, h}$ is used in the preceding sections.

One denotes by $K(h)=  ({\cal J}_h ^{-1}  \otimes I)  H(h)  (  {\cal J}_h \otimes I) $ the Hamiltonian considered as an operator in ${\cal F}_s (H_{\bf C}) \otimes {\cal H}_{sp}$.
Let us recall that $H_{ph} = h {\rm d}\Gamma (M)$. One then has,
%---
\be\label{6.19} K(h) =  ({\cal J}_h ^{-1}  \otimes I)  H(h)  (  {\cal J}_h \otimes I) =
h K_1 + h^{3/2} K_{3/2},\ee
%---
with
%---
\be\label{6.20} K_1 = {\rm d}\Gamma (M) \otimes I + I \otimes T_0,\qquad
T_0 =  \sum _{\lambda = 1}^N \sum _{m=1}^3  \beta_m  \sigma _m^{[\lambda ]}, \ee
%----6.21
\be\label{6.21} K_{3/2} = \sum _{\lambda = 1}^N \sum _{m=1}^3 \Phi_S ( a_m( x_{\lambda }) + i
b_m( x_{\lambda }) ) \otimes  \sigma _m^{[\lambda ]}
\ee
 %---
where $a_m( x_{\lambda })$ and $b_m( x_{\lambda })$ are  defined in (\ref{1.6}) and (\ref{1.7}).

The point with identifying ${\cal H}_{ph}$ with ${\cal F}_s (H_{\bf C})$ is that, the operator $ K_{1}$
and $ K_{3/2}$ are independent of $h$. Note that the operator $K_1$ maps each one of the two spaces
${\cal F }_{even } (\rho) \otimes {\cal H } _{ph}$  and ${\cal F }_{odd } (\rho) \otimes {\cal H } _{ph}$
into itself, whereas $K_{3/2}$ maps each space into the other one. This comes from the definition of the Segal field (\cite{RE-SI}).

 The following Lemma is needed for the proof of Theorem \ref{t6.5}.

\begin{lemm}\label{l6.6}  Set $u_0$
 and $\lambda _1$  defined in (\ref{6.9}).
 Then, for all  $f$ in  ${\cal F} _{odd} (\rho)\otimes {\cal H} _{sp}$,
 (resp. in ${\cal F} _{even} (\rho)\otimes {\cal H} _{sp}$), there exists
 $u$ in ${\cal F} _{odd} (\rho)\otimes {\cal H} _{sp}$,
 (resp. in ${\cal F} _{even} (\rho)\otimes {\cal H} _{sp}$) satisfying,
 %----
\be\label{6.22}( {\rm d}\Gamma (M) \otimes I   +  I \otimes (T_0  -  \lambda _1)  ) u  =  f - \Pi f,\ee
%----
 where $\Pi $ is the orthogonal projection on $u_0 $.
 In the following, we use the notation $u = (K_1 - \lambda _1)^{-1} ( f - \Pi f)$.
\end{lemm}

 Notice that $\Pi f  = 0$ if $f$ belongs to ${\cal F} _{odd} (\rho)\otimes {\cal H} _{sp}$.

{\it Proof.}  One may write,
%---
$$ f = \sum _{E \subset \{ 1, \dots , N \} } \sum _{m\geq 0}  f_{E, m }
\otimes   a_E, $$
%----
with $f_{E, m } $ in  ${\cal F} _{m} (\rho) $, vanishing for even $m$  (resp. for odd $m$
). If $ m\geq 1$, set
%----
$$  u_{E, m } (k_1 , \dots , k_m) = { f_{E, m } (k_1 , \dots , k_m) \over |k_1| + \cdots + |k_m| + 2 |\beta | |E| }.$$
%------
Since $f_{E, m } $ vanishes in neighborhood of the origin, this element is well defined, even when $E$ is empty. If $m= 0$ and $E \not = \emptyset$, set
%----
$$  u_{E, 0 }  = {  f_{E, 0 } \over  2 |\beta | |E| }. $$
%----
Then set,
 %---
$$ u =  \sum _{m\geq 1}  \sum _{E \subset \{ 1, \dots , N \} }  u_{E, m }
\otimes   a_E  +  \sum _{E \subset \{ 1, \dots , N \}\atop E \not = \emptyset  }  u_{E, 0 }
\otimes   a_E.  $$
%----
This element $u$ has the stated properties in the Lemma.

\hfill $\Box$

Theorem 6.5 follows the next Proposition.

\begin{prop}\label{p6.7} Let $\rho >0$ such that the function $\chi$ in (\ref{1.6}) and (\ref{1.7}) vanishes for $|k|\leq \rho $. One supposes $\beta \not = 0$.
 With the notations (\ref{6.4}), set $a= a_{\emptyset} \in {\cal H}_{sp}$. Then there exists
 a sequence of elements of ${\cal H}_{ph} \otimes {\cal H}_{sp} $  denoted by
 $u_{j}$ ($j\geq 0$) and a sequence of real numbers $\lambda _j$ ($j\geq 1$),  such that
 $u_0$ and $\lambda _1$ are defined in (\ref{6.9}) and for even $j$,
 %---
\be\label{6.23} u_j \in {\cal F}_{even} (\rho)  \otimes {\cal H}_{sp}\ee
 %----
 and for odd $j$,
 %---
\be\label{}  u_j \in {\cal F}_{odd } (\rho)   \otimes {\cal H}_{sp} \ee
 %----
 satisfying, for all integers $p$,
 %------
\be\label{6.25}  \Big (K(h) -(\lambda _1 h + \cdots + \lambda _{p+1} h^{p+1} ) \Big ) \ \left [ \sum _{j=0}^{2p}
 u_{j} h^{j\over 2} \right ] =\sum_{k\geq 1}   h^{ p + 1 + (k/2)} f_{2p }^{ (p + 1 + (k/2)) },\ee
 %----
 where the $f_{2p }^{ (p + 1 + (k/2)) }$ are  elements of ${\cal F }_s ( H _{\bf C} ) \otimes {\cal H}_{sp}$ and
 %----
\be\label{6.26}  (K(h) -(\lambda _1 h + \cdots + \lambda _{p+1} h^{p+1} )) \ \left [ \sum _{j=0}^{2p+1}
 u_{j} h^{j\over 2} \right ] = \sum_{k\geq 0 } h^{ p + 2 + (k/2)}
 f_{2p +1 }^{ (p + 2 + (k/2)) },\ee
 %----
where the $ f_{2p +1 }^{ (p + 2 + (k/2)) } $ are elements of  ${\cal F }_s ( H _{\bf C} )\otimes {\cal H}_{sp}$.
The sums are finite in right hand sides of  (\ref{6.25}) and (\ref{6.26}).
\end{prop}

 {\it Proof.  } In writing that the coefficient of $h^j$ ($j\leq p+1$) in the left hand side of  (\ref{6.26}) is zero, one sees that the $u_j$ and the $\lambda _j$ should satisfy the following equalities,
 %---
\be\label{6.27} ( K_1 - \lambda _1 ) u_0 = 0,\ee
 %---
\be\label{6.28} ( K_1 - \lambda _1 ) u_1 + K_{3/2} u_0 = 0,\ee
 %---
 $$ ( K_1 - \lambda _1 ) u_2 + K_{3/2} u_1 - \lambda _2 u_0 = 0, $$
 %----
 $$ ( K_1 - \lambda _1 ) u_3 + K_{3/2} u_2 - \lambda _2 u_1 = 0. $$
 %----
More generally, if $m$ is odd, one should have,
 %---
\be\label{6.29} ( K_1 - \lambda _1 ) u_m + K_{3/2} u_{m-1}  - \lambda _2 u_{m-2} - \cdots -
  \lambda _{(m+1)/2)}  u_{1} = 0\ee
  %---
and if $m$ is even,
 %---
\be\label{6.30}  ( K_1 - \lambda _1 ) u_m + K_{3/2} u_{m-1}  - \lambda _2 u_{m-2} - \cdots -
  \lambda _{(m/2) +1}  u_{0} = 0.\ee
  %---
Since
 $ {\rm d}\Gamma (M) \Psi_0 = 0$ and $(T_0 - \lambda _1) a_{\emptyset } = 0 $
 then the elements $u_0$ and $\lambda _1$ defined in (\ref{6.9}) satisfy  (\ref{6.27}).
 Since the operator $K_{3/2}$ exchanges parity then  $K_{3/2} u_0$ belongs to ${\cal F}_{odd } (\rho)   \otimes {\cal H}_{sp}$. From  Lemma 6.6, there exists $u_1$ in ${\cal F}_{odd } (\rho)   \otimes {\cal H}_{sp}$
 satisfying (\ref{6.28}). Set $p\geq 0$. Suppose that $u_0,\dots,u_{2p+1}$ and  $\lambda _1,\dots,\lambda _{p+1}$ satisfying (\ref{6.29}) and (\ref{6.30}) are constructed.
 In order to determine $u_{2p+2}$ and $\lambda _{p+2}$, one applies Lemma \ref{l6.6}
 with
 %---
\be\label{6.31}  f _{2p+2} = \left \{ \begin{matrix}  - K_{3/2} u_{2p+1} + \lambda _2 u_{2p} + \cdots + \lambda _{p+1} u_2
 & {\rm if} & p\geq 1 \\  - K_{3/2} u_{1} \hfill & {\rm if} & p= 0  \end{matrix} \right ..\ee
 %----
 Since $K_{3/2}$ exchanges parity then this element lies in  ${\cal F}_{even } (\rho)   \otimes {\cal H}_{sp}$. One defines
 $\lambda _{p+2}$ by
 %---
\be\label{6.32}   \lambda _{p+2} =- < f _{2p+2} , u_0 >.\ee
 %------
 From  Lemma \ref{l6.6}, there exists $u_{2p+2}$ in ${\cal F}_{even } (\rho)   \otimes {\cal H}_{sp}$
satisfying
  %----
$$ ( K_1 - \lambda _1 ) u_{2p+2} =  f _{2p+2} +  \lambda _{p+2} u_0, $$
%---
that is to say, (\ref{6.29}) with $p$ replaced by $p+1$. To obtain $u_{2p+3}$, one applies Lemma \ref{l6.6} with
%---
 $$ f _{2p+3} = - K_{3/2} u_{2p+2} + \lambda _2 u_{2p+1} + \cdots + \lambda _{p+2} u_1. $$
 %----
This element belongs to ${\cal F}_{odd } (\rho)   \otimes {\cal H}_{sp}$ and consequently
 $\Pi  f _{2p+3} = 0$. In view of Lemma \ref{l6.6}, there indeed exists  $u_{2p+3}$ in  ${\cal F}_{odd} (\rho)   \otimes {\cal H}_{sp}$ satisfying,
  %----
$$ ( K_1 - \lambda _1 ) u_{2p+3} =  f _{2p+3}. $$
%---
We have indeed constructed  the sequences $(u_j)$ and $(\lambda _j) $ satisfying (\ref{6.29}) and (\ref{6.30}). The
stated properties in the Proposition then follow and the proof of Theorem \ref{t6.5} is completed.

\hfill $\Box$

{\it Proof of (\ref{1.23}). }  One has,
%----
$$ <( B_m (x)\otimes I) ( u_0 + h^{1/2} u_{1} ), ( u_0 + h^{1/2} u_{1} ) >  =  2  h^{1/2}  {\rm Re}
< (B_m (x)\otimes I)   u_0 , u_{1} > .  $$
%-----
Indeed,  $B_m (x)$ maps any  element $u$ of ${\cal F} _{odd } (\rho) $  in ${\cal F} _{even } (\rho) $ and then orthogonally to $u$. The notation $\pi _{\emptyset}$ stands for the orthogonal projection in ${\cal H}_{sp}$ on the vector line spanned by $a_{\emptyset}$.
One sees,
%---
$$ < (B_m (x)\otimes I)   u_0 , u_{1} >  = < (B_m (x)\otimes I)   u_0 , (I \otimes \pi _{\emptyset} ) u_{1} >.  $$
%----
According to (\ref{6.28}),
%---
$$ (I \otimes \pi _{\emptyset} ) u_{1} =-  \Big ( {\rm d}\Gamma (M^{-1} ) \otimes \pi _{\emptyset} \Big ) K_{3/2} u_0. $$
%------
The operator ${\rm d}\Gamma (M^{-1} )$ is well defined when applied to an element of ${\cal F} _{odd } (\rho) $.
From (6.21),
%----
$$  \Big ( {\rm d}\Gamma (M^{-1} ) \otimes \pi _{\emptyset} \Big ) K_{3/2} u_0 =
 - \sum _{\lambda = 1}^N  {\rm d}\Gamma (M^{-1} ) \Phi_S ( a_3( x_{\lambda }) + i
b_3( x_{\lambda }) )\Psi_0  \otimes  a_{\emptyset}.
  $$
 %---
Therefore,
%---
$$ < (B_m (x)\otimes I)   u_0 , u_{1} >  = h^{1/2} \sum _{\lambda = 1}^N  <\Phi_S ( a_m( x) + i
b_m( x) ) \Psi_0 \  , \  {\rm d}\Gamma (M^{-1} ) \Phi_S ( a_3( x_{\lambda }) + i
b_3( x_{\lambda }) )\Psi_0  >. $$
%-------
Using standard facts concerning  Segal fields (see \cite{D-G}),
%----
$$ < \Phi_S ( a_m (x) + i b_m (x)) \Psi_0 \  ,\
{\rm d}\Gamma (M^{-1} )  \Phi_S ( ( a_3( x_{\lambda }) + i
b_3( x_{\lambda }))  )\Psi_0 > $$
%---
$$ =  {1 \over 2}\Big  ( M^{-1} a_3( x_{\lambda }) - i
 M^{-1}  b_3( x_{\lambda }) \Big )    \cdot  \Big ( a_m (x) + i b_m (x) \Big ). $$
%-----
The above scalar product is the one of  $H$.
One then deduces (\ref{1.23}). Equalities (\ref{1.24}) and (\ref{1.25}) then follows.

\hfill$\Box$

{\it Computation of $\lambda _2$.}  From (6.32) and (\ref{6.31}), one has,
%---
$$ \lambda _2 = - < f _{2} , u_0 > = < K_{3/2} u_{1} ,  u_0 > = <  u_{1} , K_{3/2} u_0 >. $$
%------
With the notations of Lemma \ref{l6.6},
%---
$$ \lambda _2 = - <  (K_1 - \lambda _1)^{-1} K_{3/2} u_0   , K_{3/2} u_0 > $$
%---
$$ = - \sum _{\lambda , \mu \leq N } \sum _{m, n \leq 3} C _{mn } ^{[\lambda , \mu]} $$
%---
$$  C _{mn } ^{[\lambda , \mu]}  =  <   (K_1 - \lambda _1)^{-1} \Phi_S ( a_m( x_{\lambda }) + i
b_m( x_{\lambda }) )  \Psi_0 \otimes \sigma _m^{[\lambda ]} a_{\emptyset }\  , \
\Phi_S ( a_n( x_{\mu }) + i
b_n( x_{\mu }) )  \Psi_0 \otimes \sigma _n^{[\mu ]} a_{\emptyset } >. $$
%-----
Suppose that $\beta = (0, 0, |\beta |)$. From the expression of the operator $(K_1 - \lambda _1)^{-1}$
(proof of Lemma \ref{l6.6}), the above bracket is zero, excepted if $m=n=3$, or if $m\leq 2 $, $n\leq 2$ and $\lambda = \mu$. In the first case, one has,
%---
$$  C _{33 } ^{[\lambda , \mu]}  = (2\pi)^{-3} \int _{\R^3}  \chi (|k|)^2  \cos ( k \cdot
(x_{\lambda } -x_{\mu } )) {k_1^2 + k_2 ^2 \over |k|^2 } dk.   $$
%------
In the second case, one gets,
%----
$$ C _{mn } ^{[\lambda , \lambda]}  =  < \sigma _m^{[\lambda ]} a_{\emptyset } , \sigma _n^{[\lambda ]} a_{\emptyset } >
(2\pi)^{-3} \int _{\R^3} { \chi (|k|)^2 |k| \over |k| + 2 |\beta |} \ { ( k \wedge e_m ) \cdot ( k \wedge e_n )
\over |k|^2 }  dk. $$
%--
The above integral is zero if $m\not = n$. Set $C = C _{1 1} ^{[\lambda , \lambda]} $.
From the preceding points,
%---
$$ \lambda _2 = - 2NC - \sum _{\lambda , \mu \leq N }   C _{33 } ^{[\lambda , \mu]}. $$
%-------

\medskip

laurent.amour@univ-reims.fr\newline
LMR EA 4535 and FR CNRS 3399, Universit\'e de Reims Champagne-Ardenne,
 Moulin de la Housse, BP 1039,
 51687 REIMS Cedex 2, France.

jean.nourrigat@univ-reims.fr\newline
LMR EA 4535 and FR CNRS 3399, Universit\'e de Reims Champagne-Ardenne,
 Moulin de la Housse, BP 1039,
 51687 REIMS Cedex 2, France.

\end{document}